\documentclass{article}

\usepackage{geometry}
\usepackage[final]{microtype}

\usepackage{import}
\usepackage[utf8]{inputenc}
\usepackage[english]{babel} 
\usepackage[T1]{fontenc} 

\usepackage{bm}
\usepackage{amssymb}
\usepackage{amsmath}
\usepackage{sansmath}
\allowdisplaybreaks
\usepackage{amsthm}
\usepackage{derivative}
\usepackage{accents}
\usepackage[scr=txupr]{mathalpha}

\newtheorem{theorem}{Theorem}
\newtheorem{corollary}{Corollary}
\newtheorem{definition}{Definition}
\newtheorem{remark}{Remark}
\newtheorem{proposition}{Proposition}

\newtheorem{problem}{Problem}

\DeclareMathOperator*{\argmin}{arg\,min}

\usepackage{enumitem}

\usepackage{lscape}
\usepackage[nice]{nicefrac}
\usepackage{mathtools}

\DeclareMathOperator{\sgn}{sgn}
\DeclareMathOperator{\spn}{span}

\DeclarePairedDelimiter\abs{\lvert}{\rvert}
\DeclarePairedDelimiter\norm{\lVert}{\rVert}

\newcommand*{\defeq}{\mathrel{\vcenter{\baselineskip0.5ex \lineskiplimit0pt
                \hbox{\scriptsize.}\hbox{\scriptsize.}}}%
                =}
                     
\def\Xint#1{\mathchoice
   {\XXint\displaystyle\textstyle{#1}}%
   {\XXint\textstyle\scriptstyle{#1}}%
   {\XXint\scriptstyle\scriptscriptstyle{#1}}%
   {\XXint\scriptscriptstyle\scriptscriptstyle{#1}}%
   \!\int}
\def\XXint#1#2#3{{\setbox0=\hbox{$#1{#2#3}{\int}$}
     \vcenter{\hbox{$#2#3$}}\kern-.5\wd0}}

\def\dashint{\Xint-}

\newcommand{\dual}[2]{\left \langle #1,#2 \right \rangle}
\newcommand{\dotp}[2]{\left ( #1,#2 \right )}

\DeclareMathAlphabet{\mathbbs}{U}{fplmbb}{m}{n}
\newcommand{\IR}{\mathbbs{R}}
\newcommand{\IC}{\mathbbs{C}}

\newcommand{\IN}{\mathbbs{N}}

\newcommand{\IT}{\mathbb{T}}
\newcommand{\ID}{\mathbb{D}}

\newcommand{\IX}{\mathbb{X}}
\newcommand{\IV}{\mathbb{V}}

\newcommand{\setT}{\mathscr{T}}
\newcommand{\setI}{\mathscr{I}}
\newcommand{\setJ}{\mathscr{J}}
\newcommand{\e}{\mathrm{e}}

\usepackage{floatrow}
\usepackage[ruled,vlined,linesnumbered]{algorithm2e}

\usepackage[bookmarks,hidelinks]{hyperref}
\usepackage[nameinlink,noabbrev]{cleveref}
\usepackage{orcidlink}
\usepackage{autonum}

\crefname{problem}{problem}{problems}
\Crefname{problem}{Problem}{Problems}
\crefname{enumi}{property}{properties}
\Crefname{Enumi}{Property}{Properties}

\usepackage{csquotes}

\numberwithin{equation}{section}

\usepackage{graphicx}
\usepackage{tikz}
\usepackage{multirow}
\usepackage{caption}
\usepackage{subcaption}
\usepackage{pgf}
\usepackage{pgfplots}
\pgfplotsset{compat=1.18}

\usepackage{tabularray}
\usepackage{booktabs}
\usepackage[square,numbers]{natbib}
\usepackage[acronym]{glossaries}

\newacronym{frb}{FRB}{Frequency Reduced-Basis}
\newacronym{pod}{POD}{Proper Orthogonal Decomposition}
\newacronym{fem}{FEM}{Finite Element Method}
\newacronym{svd}{SVD}{Singular Value Decomposition}
\newacronym[longplural={Reduced Order Models}]{rom}{ROM}{Reduced Order Model}
\newacronym[longplural={Full Order Models}]{fom}{FOM}{Full Order Model}

\newacronym[longplural={Spectral Elements}]{se}{SE}{Spectral Element}

\newacronym{rmsd}{RMSD}{Root-Mean-Square Deviation}

\newacronym{supg}{SUPG}{Streamline upwind Petrov-Galerkin}
\newacronym{vms}{VMS}{Variational Multi-Scale}

\title{The Frequency Reduced-Basis method: Reduced order models for time-dependent problems using the Laplace transform}

\author{Ricardo Reyes \orcidlink{0000-0003-0140-9564} \, \footnote{\sc{Chair of Computational Mathematics and Simulation Science, Ecole Polytechnique Fédérale de Lausanne, Switzerland}} \footnote{\sc{Department Urban Water Management, Swiss Federal Institute of Aquatic Science and Technology, Switzerland}}}

\date{}

\begin{document}
\maketitle

\begin{abstract}
In this paper, we propose the \textit{Frequency Reduced-Basis} method: a reduced-basis method to solve time-dependent partial differential equations based on the Laplace transform. Unlike traditional approaches, we begin by applying the Laplace transform to the evolution problem, yielding a time-independent boundary value problem that depends on the complex Laplace parameter. First, in an offline stage, we appropriately sample the Laplace parameter and solve the collection of problems using the finite element method.
Next, we apply a \gls{pod} to this collection of solutions in order to obtain a reduced basis that is of dimension much smaller than that of the original solution space. This reduced basis, in turn, is then used to solve the evolution problem using any suitable time-stepping method. A key insight to justify the formulation of the method relies on Hardy spaces of analytic functions. By applying the Paley-Wiener theorem we can then define an isometry between the solution of the time-dependent problem and its Laplace transform. As a consequence of this result, one may conclude that computing a \gls{pod} with samples taken in the Laplace domain produces an exponentially accurate reduced basis for the time-dependent problem.
Numerical experiments characterizing the performance of the method, in terms of accuracy and speed-up, are included for a variety of relevant time-evolution equations.
\end{abstract}



\section{Introduction}

Solving multiple-query or forecasting problems often requires the solution of multiple parametric partial differential equations quickly and reliably. Unfortunately, numerical approximations are often either computationally expensive or insufficiently accurate. At this point, the development of surrogate models that are fast and reliable becomes necessary. Reduced basis methods are suitable surrogate models, where the original high-dimensional problem is replaced by a low-rank approximation. These methods consist of two stages: an offline part, where we construct a reduced basis from a collection of `high-fidelity' solutions at selected sample points, and an online part, where we solve the low-dimensional problem defined as the projection of the original equation onto the subspace spanned by this reduced basis.

Reduced bases for time-dependent problems are usually constructed using a standard \gls{pod} \cite{Hesthaven2022,schleus_randomized_2022,baurChapterComparisonMethods2017}, which allows for a detailed data sampling. In problems where the solution can evolve without a clear pattern or at least without one that we can easily predict, we may need to sample data at a high frequency and over a long time window to be able to construct an appropriate reduced basis. If we also consider any other additional parametrization (for example physical or geometrical parameters), we may face issues sampling the data for large problems as they can become prohibitively expensive and time consuming. 

For general parametrised problems, a substantial body of work applies a Greedy sampling technique \cite{bui-thanh_model_2008,buffa_priori_2012,devore_greedy_2013}, which decreases the computational burden of sampling. Although this is ideal for most parametric problems, the evolution behaviour of a time-dependent problem does not allow us to follow a random sampling in most cases. Therefore, when we deal with parametrised time-dependent problems, it is common to find techniques that have a mix of both methods: a Greedy approach to select the sample set in the parameter space, and either an impulse response \cite{lall_subspace_2002,grepl_posteriori_2005,haasdonk_reduced_2008}, or a regular interval sampling \cite{rovas_reduced-basis_2006,nguyen_reduced_2009,grepl_certified_2012,guo_data-driven_2019,baurChapterComparisonMethods2017} in the time-domain. Although any of these approaches can give accurate results with a good enough sampling, we still face the same problem, a big computational cost for large problems. Thus, for real-world applications with computational and time constraints it is not always possible to construct an affordable and accurate reduced basis by following a traditional reduced-basis method.\footnote{Schleu\ss{} et al.\cite{schleus_randomized_2022} take an original approach to basis generation of the time-dependent problem, where rather than sampling over the global solution, the method captures a general behaviour of the problem by sampling randomly drawn starting points with random initial conditions in parallel.}\footnote{It is worth noting the work in space-time approximations where the method aims to construct a reduced basis to project the entire space-time system, rather than solving a low-dimensional system for each time step \cite{yano_space-time_2014,Choi2017}}

An alternative approach to building a reduced model consists of finding a surrogate for the transfer function in the frequency domain \cite{huynh_laplace_2011,Bigoni2020a,Bhouri2020,baurChapterComparisonMethods2017}. In this approach, mostly used in signal processing and control problems, the problem is written as an input-output linear time-invariant system, with a clearly defined transfer function. Then by using the Laplace transform we can convert the problem from the time domain to the frequency domain, find a surrogate for the transfer function and bring it back to the time domain using and inverse transform. (In this case it would be necessary to review some inverse Laplace transform methods \cite{kuhlman_review_2013}.) The major issue with this approach is how difficult and costly it can be to compute the numerical approximation of the inverse Laplace transform, which may be problem dependent and sometimes unattainable in practice.

The \textit{Frequency Reduced-Basis} method combines key ideas from frequency-domain surrogate modelling \cite{huynh_laplace_2011,peng_symplectic_2015} and standard \gls{pod}. While inspired by classical transfer-function approximations, our approach departs from them in that we avoid computing an inverse Laplace transform. Instead, we exploit the isometric structure of Hardy spaces to construct a reduced basis in the frequency domain that is then used to solve the original time-dependent problem.
A summary of the method and its main advantages is presented in \cref{ssec:Contributions}.

We intend to explain why this method works using some elements of harmonic analysis. We start by defining the approximation spaces for each problem and the equivalence between their norms: a $L^2$ Hilbert space for the time-domain problem and a $H^2$ Hardy space for the frequency-domain problem. Then we map the frequency-domain problem, set in the complex half-plane, onto the unit circle. And lastly, we justify why we only use the real part of the solution as data input using the properties of analytic functions. 

We validate the method using three numerical examples: a wave equation and a heat equation to evaluate convergence, efficiency and performance, and a advection problem to test stability.

\subsection{Method Summary} \label{ssec:Contributions}

In this work, we introduce the \gls{frb} method: a model reduction technique for time-dependent partial differential equations. The method uses the Laplace transform and the \gls{pod} to build a model reduction framework that operates in the frequency domain during the offline stage and in the time domain during the online stage.

\paragraph{Offline stage (frequency domain)}
We apply the Laplace transform to the time-dependent PDE, which yields a parametric elliptic boundary value problem in the complex Laplace domain. We then evaluate the solution at a selected set of complex frequencies and collect the corresponding snapshots. Finally, we construct a reduced basis using a \gls{svd} of these frequency-domain solutions.

\paragraph{Online stage (time domain)}
We project the original time-dependent PDE onto the reduced subspace, following the standard reduced-basis procedure. Finally, we solve the resulting low-dimensional system using a suitable time integration scheme.
\vspace{5pt}

This approach yields several advantages:
\begin{itemize}
    \item Time discretization is decoupled from basis construction, reducing sensitivity to high temporal resolution.
    \item Efficient parallelization in the offline stage due to the independence of frequency-domain solutions.
    \item Robustness for wave-like solutions, as the frequency-domain representation can capture oscillatory modes and long-time behaviour with fewer samples.
\end{itemize}

\subsection{Outline}
The paper is organized as follows. \Cref{sec:Problem} introduces the general time-dependent problem. \Cref{sec:FRB} presents the proposed Frequency Reduced-Basis method. \Cref{sec:Analysis} provides the harmonic analysis framework used in developing the method. \Cref{sec:Examples} contains numerical examples that demonstrate the method's performance. \Cref{sec:Summary} concludes with a summary of the main results.
\section{Problem model} 
\label{sec:Problem}

\subsection{Notation}
Let $\Omega \subset \mathbb{R}^d$, $d \in \{2,3\}$, be a bounded Lipschitz domain with boundary $\Gamma := \partial \Omega$.
We denote by $L^2(\Omega)$ the space of square-integrable functions on $\Omega$, equipped with the inner product
$(u, v)_{L^2(\Omega)} = \int_\Omega u(x)\,v(x)\,dx$ and the induced norm $\|u\| = \sqrt{(u,u)}_{L^2(\Omega)}$.
For $k\in\mathbb{N}$, let $H^k(\Omega)$ be the Sobolev space of functions whose weak derivatives up to order $k$ lie in $L^2(\Omega)$.
Throughout this work, we set $\mathbb{V} \defeq H^1_0(\Omega)$, the subspace of $H^1(\Omega)$ with vanishing trace on $\Gamma$, and denote its dual by $\mathbb{V}' \defeq H^{-1}(\Omega)$, with the duality pairing $\langle \cdot, \cdot \rangle_{\mathbb{V}',\mathbb{V}}$.\footnote{Note that we can set $\mathbb{V}$ as other space of functions that is adequate to the requirements of the problem, for example $H^{\textrm{div}}$, $H^{\textrm{curl}}$, or $H^1$ for different boundary conditions.}

We set the time interval $\setI=(0,t_{\textrm{f}})$ for $t_{\textrm{f}}>0$, and let $\IX$ be a separable Banach space. For each $m \in \IN_0$, we define $H^m(\setI;\IX)$ as the Bochner space of strongly measurable functions $u: \setI \rightarrow \IX$ satisfying
\begin{equation}
	\norm*{u}_{H^m(\setI ; \IX)} \defeq \left(\sum_{j=0}^m \int\limits_0^{t_{\textrm{f}}} \norm*{\partial_t^j u(t)}_{\IX}^2 \odif{t} \right)^{\frac{1}{2}} < \infty,
\end{equation}
where $\partial^j_t$ stands for the weak time derivative of order $j\in \IN_0$. In particular, if $m=0$ we set $L^2(\setI, \IX)\defeq H^0(\setI; \IX)$. We set $\IR_+ \defeq \{t \in \IR : t>0 \}$, and introduce Sobolev spaces in time $L^2(\IR_+,\IV)$ and $L^2(\IR_+,\IV')$. Given $\alpha>0$ let $L^2(\IR_+,\IV,\alpha)$ be the Hilbert space of $\IV$-valued measurable functions $u:\IR_+ \rightarrow \IV$ satisfying
\begin{align} \label{eq:normalpha}
	\norm{u}_{L^2(\IR_+,\IV,\alpha)} \defeq (u,u)^{\frac{1}{2}}_{L^2(\IR_+,\IV,\alpha)} = \left(\int\limits_{0}^{\infty} \dotp{u(t)}{u(t)}^2_{\IV} \e^{-2\alpha t} \odif{t} \right)^{\frac{1}{2}} < \infty,
\end{align}
where
\begin{align}
	\dotp{u}{v}_{L^2(\IR_+,\IV,\alpha)} \defeq \int\limits_{0}^{\infty} \dotp{u(t)}{v(t)}_{\IV} \e^{-2\alpha t} \odif{t}.
\end{align}

\subsection{Time-dependent Partial Differential Equations} \label{eq:time_dependent_pdes}

We consider a general time-dependent problem of the following form: Seek $u(\bm{x},t):\Omega\times\setI \to \IR^\ell$, where $\ell =1$ corresponds to the scalar-valued case and $\ell\in\{2,\dots,d\}$ to the vector-valued one, such that
\begin{equation} \label{eq:prob}
	\partial_{tt} u(\bm{x},t) + \mathcal{A}  u(\bm{x},t) = b(\bm{x},t), \quad \bm{x} \in \Omega,
\end{equation}
equipped with initial and boundary conditions.\footnote{Through the main body of the paper we work with the wave equation but note that the method can be used with any PDE with a defined Laplace transform.}
We define $\mathcal{A}$ as a linear operator on $\Omega$, and $\partial_{tt} u$ and $\partial_t u$ as the second and first order time derivative, respectively. Here we can set $\mathcal{A}$ as the Laplace operator $-\Delta$, an advection operator $\bm{a} \cdot \nabla$, or as any other operator. Additionally, we equip the problem with appropriate boundary conditions
\begin{equation}
	u(\bm{x},t) = u_{\Gamma}, \quad (\bm{x},t) \in \Gamma \times \setI,
\end{equation}
and suitable initial conditions
\begin{equation}
	u(\bm{x},0) = u_0(\bm{x}), \quad \text{and} \quad \partial_t u(\bm{x},0) = u_0'(\bm{x}), \quad \bm{x} \in \Omega.
\end{equation}
We also define the bilinear form $\mathsf{a}: \mathbb{V} \times \mathbb{V} \rightarrow \IR$ for $u,v \in \mathbb{V}$, for example for $\mathcal{A}=-\Delta$ as
\begin{align}
	\mathsf{a}(u,v) \defeq \int_{\Omega} \nabla u(\bm{x})^\top \nabla v(\bm{x}) \odif{\bm{x}}.
\end{align}

Now, we can introduce the variational version of the previous problem.
\begin{problem} \label{pbm:weak_problem}
Given $t \in \setI$, $u_0 \in \mathbb{V}$, and $u'_0 \in L^2(\Omega)$, find $u \in L^2(\setI,\IV)$ such that
\begin{align}
	\dual{\partial_{tt}u}{v}_{\mathbb{V}' \times \mathbb{V}} + \mathsf{a}(u,v) = \dotp{b}{v}_{L^2(\Omega)},
\end{align}
satisfying the initial conditions $u(0) = u_0$ and $\partial_t u(0) = u'_0$.
\end{problem}

\subsection{Finite Element Discretization}
\label{ssec:fe_problem}

Let $\Omega \subset \IR^d$, $d\in \{1,2,3\}$ be a bounded Lipschitz polygon with boundary $\partial \Omega$. We consider a discretization $\setT_h$ of elements in $\IR^d$ with mesh size $h$, each of them the image of a reference element $\hat{K}$ under the affine transformation $F_{K}: \hat{K} \rightarrow K$. Let us set
\begin{equation}
	\IX^{p}\left(\setT_h\right) \defeq \left\{	w \in \mathcal{C}^0(\Omega) \mid \forall K \in \setT_h: \left. w\right|_K \circ F_K \text { is a polynomial of degree } p \right\}.
\end{equation}
We set $\IV_h = \IV \cap \IX^p(\setT_h)$ as a conforming finite-dimensional subspace of $\IV$, with $n_h = \dim(\IX^{p} \left(\setT_h\right))$ and a basis $\{\varphi_1,\dots,\varphi_{n_h}\}$. We define $\mathsf{P}_h: \IV \rightarrow \IV_h$ to be the projection operator into $\IV_h$. To write the discrete problem, let us consider the solution ansatz
\begin{equation}
	u_h(\bm{x},t) = \sum_{j=1}^{n_h} \varphi_j(\bm{x}) \mathsf{u}_j(t), \quad\bm{x} \in \Omega, \quad t\in (0,t_{n_t}].
\end{equation}

\begin{problem}[Semi-discrete problem] \label{pr:sdp}
Let $u^0, u^1 \in \IV$, and $b \in \mathcal{C}^0(\setI ; L^2(\Omega))$. We seek $u_h \in \mathcal{C}^2(\setI, \IV_h)$ such that for each $t \in \setI$ it holds
\begin{align}\label{eq:semi_discrete}
	\odv[order=2]{}{t} \dotp{u_h}{v_h}_{L^2(\Omega)} + \mathsf{a} \left(u_h,v_h\right) = \dotp{b}{v_h}_{L^2(\Omega)}, \quad \forall t \in (0,t_{n_t}], \quad \text{and} \quad \forall v_h \in \IV_h,
\end{align}
with initial conditions $u_h(0) = \mathsf{P}_{h} u^0\in \IV_h$ and $ \odv{}{t} u_h(0) = \mathsf{P}_{h} u^1\in \IV_h$.
\end{problem}

To describe the fully discrete version of \cref{pr:sdp}, we consider the discrete time interval $\setI_h=\{t_1,\ldots,t_{n_t} \}$ with $n_t$ total time steps. Let us set $\bm{\mathsf{u}}(t_n) = (\mathsf{u}_1(t_n),\dots, \mathsf{u}_{n_h}(t_n))^{\top} \in \IR^{n_h}$, and likewise $\ddot{\bm{\mathsf{u}}}(t_n)$ and $\dot{\bm{\mathsf{u}}}(t_n) \in \IR^{n_h}$ as the approximations of $\partial_{tt} u$ and $\partial_{t} u$, respectively. We define $\bm{\mathsf{M}} \in \IR^{n_h \times n_h}$ and $\bm{\mathsf{A}} \in \IR^{n_h \times n_h}$ as
\begin{equation}
	(\bm{\mathsf{M}})_{i,j} = \dotp{\xi_i}{\xi_j}_{L^2(\Omega)} \quad \text{and} \quad (\bm{\mathsf{A}})_{i,j} = \mathsf{a} \left(\xi_i,\xi_j \right), \quad i,j \in \{1,\dots,n_h\},
\end{equation}
referred to as the mass matrix and the stiffness matrix of the bilinear form $\mathsf{a}(\cdot,\cdot)$.

As a time integration algorithm, we follow the generalized $\alpha$-method proposed in \cite{chungTimeIntegrationAlgorithm1993}.
\begin{problem}[Fully-discrete problem] \label{pr:dp}
We seek $\bm{\mathsf{u}}_{n+1} = \bm{\mathsf{u}}(t_{n+1}) \in \IR^{n_h}$, such that
\begin{equation}
    \bm{\mathsf{M}} \ddot{\bm{\mathsf{u}}}_{n+1-\alpha_{m}} + \bm{\mathsf{A}} \bm{\mathsf{u}}_{n+1-\alpha_{f}} = \bm{\mathsf{b}}_{n+1-\alpha_{f}},
\end{equation}
with the notation $\bm{\mathsf{b}}_{n+1-\alpha} = (1-\alpha) \bm{\mathsf{b}}_{n+1} + \alpha \bm{\mathsf{b}}_{n}$, by using the approximation
\begin{align}
    \bm{\mathsf{u}}_{n+1} &= \bm{\mathsf{u}}_n + \Delta t \dot{\bm{\mathsf{u}}}_n + \frac{\Delta t^2}{2} \left( (1-2\beta) \ddot{\bm{\mathsf{u}}}_n + 2 \beta \ddot{\bm{\mathsf{u}}}_{n+1} \right) \\
    \dot{\bm{\mathsf{u}}}_{n+1} &= \dot{\bm{\mathsf{u}}}_{n} + \Delta t \left( (1-\gamma) \ddot{\bm{\mathsf{u}}}_n +\gamma \ddot{\bm{\mathsf{u}}}_{n+1} \right).
\end{align}
We set the parameters as $\alpha_m, \alpha_f \leq \frac{1}{2}$ and $\gamma = \frac{1}{2} - \alpha_m + \alpha_f$, $\beta \geq \frac{1}{4} + \frac{1}{2} (\alpha_f-\alpha_m)$ which ensures that the method is unconditional stable and of second-order accuracy.
\end{problem}

\subsection{Construction of the reduced basis: Time-dependent approach}

The traditional approach to model reduction for time-dependent problems consists in finding a reduced space $\IV_h^{\textrm{R}}\subset\IV_h$ with basis $\{\varphi_1^\textrm{R},\dots,\varphi_{n_r}^\textrm{R}\}$, where the dimension of the reduced space $n_r$ is considerably smaller compared to that of the `high fidelity' space $\IV_h$. We are interested in construction of a reduced basis for the discrete solution manifold
\begin{equation} \label{eq:tsnap}
    \bm{\mathsf{S}} = \{ \bm{\mathsf{u}}(t_j) \; | t_j \in \setI_h \} \in \IR^{n_h \times n_t}.
\end{equation}
As mentioned in the introduction, the most common approach to constructing the reduced basis for a time-dependent problem consists of following a \gls{pod} to find a basis $\Phi = \{\phi_1,\ldots,\phi_{n_r}\}$ such that
\begin{equation} \label{eq:tPOD}
	\Phi = \min_{\substack{\Phi \in \IR^{n_h \times n_r}: \\ \Phi^{\top} \Phi = \bm{I}_{n_r \times n_r}}}
    	\sum_{j=1}^{n_t} \norm*{\bm{\mathsf{u}}(t_j) - \sum_{k=1}^{n_r} \dotp{\Phi_k}{\bm{\mathsf{u}}(t_j)} \Phi_k
	}^2_{\IR^{n_h}}
\end{equation}

Obtaining these sampled solutions requires solving $n_t$ linear systems of size $n_h$, which becomes increasingly expensive for refined meshes, long time intervals, or high-frequency regimes. Since this work focuses on the time-dependent component of the problem, we do not address techniques aimed at reducing the computational cost associated with spatial discretization or parameter sampling. However, it is important to mention that to reduce these computational costs we can follow works in greedy algorithms \cite{grepl_posteriori_2005,Rozza2008a} for the parameter space, and works in adaptive mesh refinement \cite{Carlberg2014,Etter2019} and incomplete sampling \cite{everson_karhunenloeve_1995,venturi_gappy_2004,willcox_unsteady_2006,Peherstorfer2016} for the spatial discretization.

\section{The Frequency Reduced-Basis method}
\label{sec:FRB}

In this section we described an alternative method for the construction of the reduced basis. Instead of directly performing a time discretization of \cref{pr:sdp}, we use the Laplace transform $\mathcal{L}$ to describe the problem in the frequency domain. We can define the Laplace transform of a function $f$ as
\begin{equation} \label{eq:laplace}
	\IC_+ \ni s \mapsto \widehat{f}(s) \defeq \langle f(t),\e^{-st}\rangle = \int_{0}^{\infty}\e^{-st}f(t) \odif{t}.
\end{equation}
where $\IC_+ \defeq \{ s \in \IC: \Re(s) >0 \}$ denotes the positive real-part half of the complex plane. Using the properties of the Laplace transform, we can define the composition with a steady state operator as $\mathcal{L}\{\mathcal{A}f (t)\} = \mathcal{A}\widehat{f}(s)$, and the Laplace transform of the $n$-order derivative $f^{(n)}$ as 
\begin{equation}
    \mathcal{L}\{f^{(n)}\}=s^n\widehat{f}(s)-\sum_{k=1}^n s^{n-k}f^{(k-1)}(0).
\end{equation}

Assuming that the source term satisfies the space-time separation ---$b(\bm{x},t)=b_x(\bm{x})b_t(t)$, that we have properly defined the Hardy space $H^2$ ---which we discuss in more detail in the next section, and that the Laplace transform of the $\bm{\mathsf{u}}$ is applied component-wise ---$\widehat{\bm{\mathsf{u}}}(s) = \mathcal{L}\left\{ {\bm{\mathsf{u}}}(t)\right\}$; we can write a frequency-domain version of \cref{pr:dp}.
\begin{problem}[Frequency discrete problem] \label{pr:fdp}
Seek $\widehat{\bm{\mathsf{u}}}_h \in H^2(\IC_+,\IV_h)$ such that
\begin{equation}
    \left(s^2 \bm{\mathsf{M}} + \bm{\mathsf{A}} \right) \widehat{\bm{\mathsf{u}}}(s) = \widehat{\bm{\mathsf{b}}}(s) - s \bm{\mathsf{u}}_{0,h} - \bm{\mathsf{u}}'_{0,h}.
\end{equation}
\end{problem}
Here, it is important to note that the spaces $\IV$ and $\IV'$ for the frequency domain solutions are complex valued. However we use the same notation as for the real-valued counterpart for simplicity. 

Analogously to the naïve construction of the reduced basis using time solutions in \cref{eq:tPOD}, we can construct the reduced basis following a \gls{pod}, but in this case over frequency domain. 
To this end we sample over a set of complex frequencies $\setJ_h = \left\{s_1,\dots,s_{n_s}\right\} \subset \IC_{+}$ with $n_s$ total frequencies, to build the manifold
\begin{equation} \label{eq:fsnap}
    \widehat{\bm{\mathsf{S}}} = \{ \widehat{\bm{\mathsf{u}}}(s_j) \; | s_j \in \setJ_h \} \in \IC^{n_h \times n_s},
\end{equation}
from which we can construct a basis $\Phi = \{\phi_1,\ldots,\phi_{n_r}\}$ such that
\begin{equation} \label{eq:fPOD}
	\Phi = \min_{\substack{\Phi \in \IR^{n_h \times n_r}: \\ \Phi^{\top} \Phi = \bm{I}_{n_r \times n_r}}}
    	\sum_{j=1}^{n_s} \norm*{\widehat{\bm{\mathsf{u}}}(s_j) - \sum_{k=1}^{n_r} \dotp{\Phi_k}{\widehat{\bm{\mathsf{u}}}(s_j)} \Phi_k
	}^2_{H^2(\IC_+,\IV_h)}
\end{equation}

Obtaining this frequency solution manifold requires the computation of $n_s$ linear systems of size $n_h$. There are two key differences with the time-dependent approach: the computation of the linear systems can be done in parallel as they do not need to be solved sequentially, and we may need fewer samples in the frequency-domain to construct the reduced basis $(n_s < n_t)$ as in many cases a complex temporal solution can be represented with a small set of frequencies. We solve the minimization problem in \cref{eq:fPOD} by computing the singular value decomposition $\widehat{\bm{\mathsf{S}}} = \sum_k^{n_r} \sigma_k \Phi_k \bm{\psi}_k^*$, and by using the first $n_r$ left-singular vectors we can define the reduced basis vectors as
\begin{equation}
	\varphi^{\textrm{R}}_k = \sum_{j=1}^{n_h} \left( \Phi_k	\right)_j \varphi_{j},
\end{equation}
and the reduced space as $\mathbb{V}^{\textrm{R}}_{h} = \spn \left\{ \varphi^{\textrm{R}}_1, \dots, \varphi^{\textrm{R}}_{n_r} \right\}$.

So far, we have described the basis construction method, but we still have not defined the Hardy space $H^2$, or the sampling set $\setJ_h$. For now let us define the sampling set as $$\setJ_h = \{ s_j =\alpha + \iota \beta \cot{\theta_j} | \theta_j \in (0,\pi) \}$$ with $\alpha, \beta \in \IR_+$. In \cref{sec:Analysis} we provide these definitions in detail and explain how we developed the method. \Cref{alg:basis} summarizes the proposed frequency reduced basis method.
\begin{algorithm}
\caption{Construction of the reduced basis using the \textit{Frequency Reduced-Basis} method} \label{alg:basis}
    Compute the Laplace $\mathcal{L}$ transform of the original time-dependent problem
    
    \For{$s_n =\alpha + \iota \beta \cot{\theta_n} \in \setJ_h | \theta_n \in (0,\pi)$}
        {Solve the frequency-domain problem for $\widehat{u}(\bm{x},s)$}
        
    Set the complex-valued data set. $\widehat{\bm{\mathsf{S}}} = \{ w_1\widehat{\bm{\mathsf{u}}}(s_1), \dots, w_{n_s}\widehat{\bm{\mathsf{u}}}(s_{n_s}) \} \in \IC^{n_h \times n_s}$, with $w_j = \frac{2\beta}{\sin^2 \theta_j}$ 

    Construct the basis $\Phi$ using a \gls{svd}. $\widehat{\bm{\mathsf{S}}} = \Phi \Sigma \Psi^*$
\end{algorithm}

\subsection{Time-dependent reduced order model} \label{ssec:rom}

After constructing the reduced basis, we project \cref{pr:sdp} onto the reduced space $\IV_h^{\textrm{R}}$, following the standard model reduction procedure. Let $\mathsf{P}_h^{\textrm{R}}: \IV \rightarrow \IV_h^{\textrm{R}}$ denote the projection operator. The reduced solution is then given by the ansatz
\begin{equation}
    u_h^{\textrm{R}}(\bm{x},t) = \sum_{j=1}^{n_h} \varphi_{j} \sum_{k=1}^{n_r} \phi_{j,k} \mathsf{u}_k(t) = \sum_{j=1}^{n_h}  \sum_{k=1}^{n_r} \varphi_{j,k}^{\textrm{R}} (\bm{x}) \mathsf{u}_k(t) 
\end{equation}

\begin{problem}[Reduced semi-discrete problem] \label{pr:sdpr}
We seek $u^{\textrm{R}}_{h} \in \mathcal{C}^2(\setI, \IV_h^{\textrm{R}})$ such that for each $t \in \setI$ it holds
\begin{align}\label{eq:semi_discrete_rom}
	\frac{d^2}{dt^2} \dotp{u_h^{\textrm{R}}}{v_h^{\textrm{R}}}_{L^2(\Omega)} + \mathsf{a} \left(u_h^{\textrm{R}},v_h^{\textrm{R}} \right) = \dotp{b}{v_h^{\textrm{R}} }_{L^2(\Omega)}, \quad \forall t \in (0,t_{n_t}], \quad \text{and} \quad \forall v_h^{\textrm{R}} \in \IV_{h}^{\textrm{R}}
\end{align}
with initial conditions $u_h(0) = \mathsf{P}_{h}^{\textrm{R}} u^0\in \IV_h^{\textrm{R}}$ and $ \odv{}{t} u_h(0) = \mathsf{P}_{h}^{\textrm{R}} u^1\in \IV_h^{\textrm{R}}$.
\end{problem}

In the same way as for the high fidelity solution in \cref{pr:dp}, we can write a fully-discrete problem for the reduced order model using the same time integration scheme. 
We set $\bm{\mathsf{u}}^{\textrm{R}}(t) = (\mathsf{u}^{\textrm{R}}_1(t),\dots, \mathsf{u}^{\textrm{R}}_{n_r}(t))^{\top} \in \IR^{n_r}$, and likewise $\ddot{\bm{\mathsf{u}}}^{\textrm{R}}(t)$ and $\dot{\bm{\mathsf{u}}}^{\textrm{R}}(t) \in \IR^{n_r}$. We also define the reduced mass matrix $\bm{\mathsf{M}}^{\textrm{R}} \in \IR^{n_r \times n_r}$, the reduced stiffness matrix $\bm{\mathsf{A}}^{\textrm{R}} \in \IR^{n_r \times n_r}$, and the reduced source term as
\begin{equation}
	\bm{\mathsf{M}}^{\textrm{R}} = \Phi^\top \bm{\mathsf{M}} \Phi, \quad \bm{\mathsf{A}}^{\textrm{R}} = \Phi^\top \bm{\mathsf{A}} \Phi, \quad \text{and} \quad \bm{\mathsf{b}}^{\textrm{R}} = \Phi^\top \bm{\mathsf{b}}
\end{equation}
with $\Phi \in \IR^{n_h \times n_r}$ the discrete reduced basis.

\begin{problem}[Reduced fully-discrete problem] \label{pr:dpr}
We seek $\bm{\mathsf{u}}^{\textrm{R}}_{n+1} = \bm{\mathsf{u}}^{\textrm{R}}(t_{n+1}) \in \IR^{n_r}$, such that
\begin{equation}
    \bm{\mathsf{M}}^{\textrm{R}} \ddot{\bm{\mathsf{u}}}^{\textrm{R}}_{n+1-\alpha_{m}} + \bm{\mathsf{A}}^{\textrm{R}} \bm{\mathsf{u}}^{\textrm{R}}_{n+1-\alpha_{f}} = \bm{\mathsf{b}}^{\textrm{R}}_{n+1-\alpha_{f}},
\end{equation}
using the same approximations and parameters as for the high fidelity problem.
\end{problem}

\section{Elements of harmonic analysis} \label{sec:Analysis}

We recall the most important elements of harmonic analysis that are relevant for the subsequent analysis. We follow \cite{partingtonBanachSpacesAnalytic,partingtonLinearOperatorsLinear2004,ricciHardySpacesOne}, \cite[Ch. 3, 7, 8]{hoffman_banach_2014}, \cite[Ch. 4]{RR97} and \cite[Section 6.4]{hille1996functional}. The main goal of this section is to state basic definitions, notation, and further properties of Hardy classes of vector-valued holomorphic functions, i.e. with values in either Hilbert or Banach spaces.

\subsection{The Laplace Transform and the Hardy space}

First, let us introduce an appropriate space for $\widehat{f}$. We set
\begin{align}
	\IC_\alpha \defeq\{s \in \mathbb{C}: \,\Re(s) > \alpha \}.
\end{align}
where $\alpha$ ensures that the singularities of $\widehat{f}$ are to the left of the semiplane $\IC_\alpha$.
\begin{definition}[Hardy space in $\IC_\alpha$] \label{df:Hardy-plane}
    We denote $H^p(\IC_\alpha,\IV)$ the space of holomorphic functions $\widehat{f}: \Omega \rightarrow \IV$ in the positive half complex plane. For $1\leq p \leq \infty$ and any Banach space $\IV$ we define the norm
    \begin{equation} \label{eq:norm-c}
        \norm{\widehat{f}}_{H^p(\IC_\alpha, \IV)} = \sup_{\sigma>\alpha} \left(\int_{-\infty}^{\infty} \norm{\widehat{f}(\sigma + \imath \tau)}_{\IV}^p \odif{\tau} \right)^{\nicefrac{1}{p}} < \infty
    \end{equation}
\end{definition}
\begin{proposition}\label{prop:hardy-plane-limit1}
Let $p\in [1,\infty)$ and $\alpha>0$. For $f \in H^p(\IC_\alpha,\IV)$, the function 
	\begin{equation}
		T(\sigma,f) = \int\limits_{-\infty}^{\infty} \norm*{\widehat{f}(\sigma+\imath \tau)}^2_\IV \odif{\tau}
	\end{equation}
	is continuous monotone decreasing function of $\sigma$ for $\sigma \geq \alpha$. 
	In particular, $T(\alpha,f) = \norm{f}^p_{H^p(\IC_\alpha,\IV)}$ and $\displaystyle\lim_{\sigma \rightarrow \infty} T(\sigma,f)  = 0$.
\end{proposition}

\begin{proposition} \label{prop:hardy-plane-limit2}
    The function $\widehat{f}(\sigma+\imath \tau)$ has boundary values for $\widehat{f}_{\alpha}(\imath \tau) = \lim_{\sigma\rightarrow \alpha^{+}} \widehat{f}(\sigma+\imath \tau)$ almost everywhere. The boundary function $\widehat{f}_{\alpha}$ lies in $L^p(\imath\IR,\IV)$ and satisfies $\lVert \widehat{f}_{\alpha} \rVert_{L^{p}(\imath\IR,\IV)} = \lVert \widehat{f} \rVert_{H^{p}(\IC_\alpha,\IV)}$. By identifying $\widehat{f}$ with $\widehat{f}_{\alpha}(\imath \tau)$, $H^p(\IC_\alpha,\IV)$ can be regarded as a closed subspace of $L^p(\imath\IR,\IV)$. 
\end{proposition}

\subsection{Isometry between time and frequency domains}

\begin{theorem} [Paley-Wiener] \label{th:PW}
    $\widehat{f} \in L^2(\imath\IR)$ if and only if there exists a function $f \in L^2(\IR)$ such that $\widehat{f}(s) = \mathcal{L}(f)$.
    We define $L^2(\IR)=L^2(-\infty,0)\oplus L^2(\IR_+)$ by decomposing the function into its values $t<0$ and $t>0$. Then, by this decomposition we see
    \begin{equation}
       (L^2(\IR)=L^2(\mathbb{R}_-)\oplus L^2(\IR_+)) \xrightarrow{\makebox[0.5cm]{$\mathcal{L}$}} (L^2(i\IR)=H^2(\IC_-)\oplus H^2(\IC_+)). 
    \end{equation}
    We can define the Laplace transform as the map $\mathcal{L}:L^2(\IR_+;\IV) \mapsto H^2(\IC_+;\IV)$.
\end{theorem}

\begin{corollary}\label{cor:paley_wiener_alpha}
\begin{equation}
	\norm{f}_{L^2(\mathbb{R}_+,\IV,\alpha)} = \frac{1}{\sqrt{2\pi}} \norm*{\mathcal{L}\{f\}}_{H^2(\IC_\alpha,\IV)}
\end{equation}
\end{corollary}

\begin{proof}
For $f \in L^2(\IR_+,\IV,\alpha)$ with $\alpha>0$ we set $f_\alpha(t) = \e^{-\alpha t}f(t)$. Recalling the properties of the Laplace transform
we have that
\begin{equation}\label{eq:shift_laplace_transform}
	\mathcal{L} \{f_\alpha\}(s)	= \mathcal{L}\{f\}(s+\alpha), \quad \Re\{s\} > 0.
\end{equation}
We observe that $f_\alpha \in L^2(\IR_+,\IV)$. Hence, it follows from \cref{th:PW} that
\begin{equation}\label{eq:paley_wiener_alpha}
	\norm{f_\alpha}_{L^2(\mathbb{R}_+, \IV)} = \frac{1}{\sqrt{2\pi}} \norm*{\mathcal{L}\{f_\alpha\}}_{H^2(\IC_+,\IV)}.
\end{equation}
By shifting the vertical integration line in the definition of
$H^2(\IC_\alpha,V)$ and recalling \eqref{eq:shift_laplace_transform}, we obtain 
\begin{equation}\label{eq:shift_frequency}
\begin{aligned}
	\norm{ \mathcal{L}\{f_\alpha\}}_{H^2(\IC_+,\IV)} &= \sup_{\sigma>0} \left(\int_{-\infty}^{-\infty} \norm{\mathcal{L}\{f_\alpha\}(\sigma+\imath \tau)}^2_\IV  \odif{\tau} \right)^{\frac{1}{2}} \\
    &= \sup_{\sigma>0} \left( \int_{-\infty}^{-\infty} \norm{\mathcal{L} \{f \}(\sigma+\alpha+\imath\tau) }^2_\IV \odif{\tau} \right)^{\frac{1}{2}} \\
	&= \sup_{\sigma>\alpha} \left( \int_{-\infty}^{-\infty} \norm{\mathcal{L} \{f\}(\sigma+\imath\tau)}^2_\IV \odif{\tau} \right)^{\frac{1}{2}} \\
	&= \norm{ \mathcal{L}\{f\} }_{H^2(\IC_\alpha,\IV)}.
\end{aligned}
\end{equation}
From \cref{eq:normalpha}, we also know that 
\begin{equation}\label{eq:shift_time}
	\norm{f_\alpha}_{L^2(\mathbb{R}_+,\IV)} = \norm{f}_{L^2(\mathbb{R}_+,\IV,\alpha)}.
\end{equation}
By combining \cref{eq:shift_time} and \cref{eq:shift_frequency}, together with \cref{eq:paley_wiener_alpha}, one obtains the final result.
\end{proof}

\subsection{Hardy spaces in the unit disk \texorpdfstring{$\ID$}{D} and the unit circle \texorpdfstring{$\IT$}{T}}

Since the behaviour of $\widehat{f}_\alpha$ can be highly oscillatory as
$\Im \{s\}\rightarrow\infty$ and therefore \cref{eq:norm-c} would be difficult to compute, we apply a transformation to map $s$ to a new complex variable that has a better behaviour. Let $\ID = \{ \omega \in  \IC :\norm{\omega}<1 \}$ denote the unit disk in the complex plane and $\IT = \{ \omega \in \IC :\norm{\omega}=1 = \e^{i\theta}: 0\leq \theta \leq 2\pi \}$ the boundary of $\ID$. There is a natural isomorphism between the Hardy spaces in $\IC_+$ and $\ID$, and by extension between $\IC_\alpha$ and $\ID$. We apply the Möbius transformation $\mathcal{M}:\IC_\alpha \mapsto \ID$, given by $\mathcal{M}(s) = \omega = \frac{s-\alpha-\beta}{s-\alpha+\beta}$ to map $\widehat{f}_\alpha(s)$ in the half plane to the unit disk.

\begin{figure}[ht]
    \centering
    \includegraphics[scale=0.9]{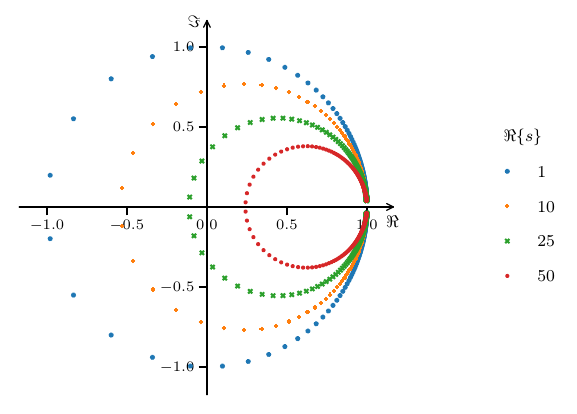}
    \caption{Illustration of the Möbius transformation with $\alpha=1$ and $\beta=30$.}
    \label{fig:mobius}
\end{figure}

By mapping the domain from $\IC_\alpha$, we ensure that the singularities of $\widehat{f}(s)$ in the half-plane are mapped to the exterior of the unit disk $\ID$. An exception occurs when the singularity is located at infinity. By increasing $\beta$ we extend the sampling to higher frequencies (approaching $\infty$ in the Bromwich integral) and therefore we should carefully choose this parameter in function of the mesh size and other parameters that we discuss later on.
\Cref{fig:mobius} shows an illustration of the mapping from different contours $\Re s \geq \sigma$, and how the contours approach to $1$ when $\Re s \rightarrow \infty$. 

Let us now define the appropriate space for the function $\widehat{f}$ in the unit disk and some theorems regarding the norm and the interoperability between the Hardy space in the disk and the half plane.
\begin{definition}[Hardy space in $\ID$] \label{df:Hardy-disk}
     We denote $H^p(\ID,\IV)$ the space of holomorphic functions $\widehat{g}: \Omega \rightarrow \IV$ in the unit disc for which the functions $\widehat{g}(r\e^{\imath \theta})$ are bounded in $L^p$-norm as $r \rightarrow 1$. For $1\leq p \leq \infty$ and any Banach space $\IV$ we define the norm
    \begin{equation} \label{eq:HC-norm}
        \norm{\widehat{g}}_{H^p(\ID, \IV)} = \sup_{0 \leq r < 1} \left( \frac{1}{2\pi} \int_{0}^{2\pi} \norm{\widehat{g}(r\e^{i\theta})}_{\IV}^p \odif{\theta} \right)^{\nicefrac{1}{p}} < \infty
    \end{equation}
\end{definition}
\begin{theorem}[Fatou's theorem] \label{th:Fatou}
    The functions $\widehat{g} \in H^p(\ID)$ have a radial limit $$\widehat{g}_\circ (\e^{i\theta}) = \lim_{r\rightarrow 1} \widehat{g}(r\e^{i\theta})$$ that exist almost everywhere in $\IT$, with $\widehat{g}_\circ \in L^p(\IT,\IV)$, and $\norm{\widehat{g}}_{H^p(\IC_\alpha,\IV)} = \norm{\widehat{g}_\circ}_{L^p(\IT,\IV)}$. We also identify $\widehat{g}$ with $\widehat{g}_\circ$, setting $H^p(\ID,\IV)$ as a closed subspace of $L^p(\IT,\IV)$ and a Banach space. A proof of this theorem is available in \cite{RR97}
\end{theorem}

\begin{theorem} \label{th:Hoffman}
    A function $\widehat{g}$ defined on $\IT$ is in $L^p(\IT,\IV)$ if and only if the function $\widehat{f}_{\alpha}: \imath \IR \mapsto \IC$ defined by 
    \begin{equation} \label{eq:mapDH}
        \widehat{f}_{\alpha}(s) = \frac{\beta^{\nicefrac{1}{p}}}{\pi^{\nicefrac{1}{p}} (\beta-\alpha+s)^{\nicefrac{2}{p}}}\widehat{g}(\mathcal{M}(s))
    \end{equation}
    is in $L^p(\imath\IR,\IV)$. Setting $\mathcal{M}^{-1}(\omega) = s = \frac{2\beta}{1-\omega}+\alpha -\beta$ we also verify
    \begin{equation} \label{eq:mapHD}
        \widehat{g}(\omega) = \frac{(4 \beta \pi)^{\nicefrac{1}{p}}}{(1-\omega)^{\nicefrac{2}{p}}}\widehat{f}_{\alpha}(\mathcal{M}^{-1}(\omega)),
    \end{equation}
    and moreover $\norm{\widehat{f}_{\alpha}}_{L^p(\imath\IR,\IV)} = \norm{\widehat{g}}_{L^p(\IT,\IV)}$.
\end{theorem}
\begin{proof}
    We follow a similar proof as the one in \cite{partingtonBanachSpacesAnalytic}, for a similar Möbius transformation. First we check that $\abs{\omega} <1$ if and only if $\Re \mathcal{M}^{-1}(s) >0$, and that $\mathcal{M}^{-1} \circ \mathcal{M} = I$. Then we have
    \begin{align}
        \norm*{\widehat{f}_{\alpha}}_{L^p(\imath\IR, \IV)}^p = \int_{-\infty}^\infty \norm*{\widehat{f}(s)}^p_{\IV} \odif{s} = - \imath \int_{\imath\IR} \frac{\beta}{\pi \abs{\beta-\alpha+s}^2} \norm*{\widehat{g}(\mathcal{M}(s))}^p_{\IV} \odif{s}.
    \end{align}
    Using $\mathcal{M}^{-1}(s)$, we replace $\odv{s}{\omega}=-\frac{2\beta}{(1-\omega)^2}$, and $\beta-\alpha+s = \frac{2\beta}{1-\omega}$. Noting that $\frac{\abs{1-\omega}^2}{(1-\omega)^2} = -\frac{1}{\omega}$ for $\omega \in \IT$, we get
    \begin{align}
        \norm*{\widehat{f}_{\alpha}}_{L^p(\imath\IR, \IV)}^p = - \imath \int_{\abs{\omega}=1} \frac{\abs{1-\omega}^2 \beta}{4\beta^2 \pi} \norm*{\widehat{g}(\omega)}^p_{\IV} \frac{2\beta}{(1-\omega)^2} \odif{s} = - \frac{1}{2\pi} \imath \int_{\abs{\omega}=1} \frac{1}{\omega} \norm*{\widehat{g}(\omega)}^p_{\IV} \odif{\omega}.
    \end{align}
    Lastly, writing $\omega=\e^{\imath\theta}$ and $\odif{\omega} = \imath \e^{\imath \theta} d\theta$, we get
    \begin{equation}
        \norm*{\widehat{f}_{\alpha}}_{L^p(\imath\IR, \IV)}^p = - \frac{1}{2\pi} \imath \int_{0}^{2\pi} \frac{\imath \e^{\imath \theta}}{\e^{\imath\theta}} \norm*{\widehat{g}(\e^{\imath \theta})}^p_{\IV} \odif{\theta} = \frac{1}{2\pi} \int_{0}^{2\pi} \norm*{\widehat{g}(\e^{\imath \theta})}^p_{\IV} \odif{\theta}.
    \end{equation}
\end{proof}

\begin{corollary} \label{cor:isometry}
Let $\widehat{f} \in H^2(\IC_\alpha)$. The mapping $\mathcal{D}: H^2(\IC_\alpha) \mapsto H^2(\ID)$ defined by 
\begin{equation} \label{eq:isometry}
    \mathcal{D}\left( \widehat{f} \right)(\e^{\imath \theta}) =  \frac{2 \sqrt{\beta \pi}}{1-\e^{\imath \theta}}\widehat{f} \left(\mathcal{M}^{-1}(\e^{\imath \theta}) \right)
\end{equation}
is an isometric isomorphism.
\end{corollary}

\begin{proof}
    A direct consequence of \cref{prop:hardy-plane-limit2,th:Fatou} is the isometry between the Hardy spaces in the disk and the half-plane and its boundary values. Moreover, \cref{th:Hoffman} is easy to proof for $p=2$ by replacing \cref{eq:isometry} and $\mathcal{M}^{-1}(\e^{\imath \theta})$ in \cref{eq:mapHD} as done for \cref{eq:mapDH}.
    We also compute the norm of the map $\mathcal{D}$ in \cref{eq:isometry} as
    \begin{align}
        \norm*{\mathcal{D} \left( \widehat{f} \right)}_{L^2(\IT,\IV)} &= \frac{1}{2\pi} \int_0^{2\pi} \frac{4 \beta \pi}{\abs*{1-\e^{\imath \theta}}^2} \norm*{\widehat{f}(\mathcal{M}^{-1}(\e ^{\imath\theta}))}^2_{\IV} \odif{\theta} \\
        &= \int_0^{2\pi} \frac{2 \beta}{\abs*{1-\e^{\imath \theta}}^2} \norm*{\widehat{f} \left( \alpha + \beta \frac{1+\e^{\imath \theta}}{1-\e^{\imath \theta}} \right)}^2_{\IV} \odif{\theta} \\
        &= \int_0^{\pi} \frac{2\beta}{\sin^2 \theta} \norm*{\widehat{f} \left( \alpha + \beta \cot \theta \right)}^2_{\IV} \odif{\theta}  \label{eq:int-rule}
\end{align}
\end{proof}

\subsection{The Frequency Reduced-Basis method} \label{ssec:FRB-harmonic}

Using the isometries presented in \cref{th:PW,cor:paley_wiener_alpha} we justify the use of the \gls{pod} in the frequency domain instead of the standard time-domain one in \cref{eq:tPOD}.
Then, using \cref{eq:int-rule} we write a \gls{pod} analogous to \cref{eq:fPOD} but this time evaluated in functions over the unit circle $\IT$. We want to find the reduced basis $\Phi = \{\phi_1,\ldots,\phi_{n_r}\}$ such that
\begin{equation} \label{eq:fPOD-circle}
	\Phi = \argmin_{\substack{\Phi \in \IR^{n_h \times n_r}: \\ \Phi^{\top} \Phi = \bm{I}_{n_r \times n_r}}}
    	\sum_{j=1}^{n_s} w_j \norm*{\widehat{\bm{\mathsf{u}}}(s_j) - \sum_{k=1}^{n_r} \dotp{\phi_k}{\widehat{\bm{\mathsf{u}}}(s_j)} \phi_k
	}^2_{L^2(\IT,\IV_h)},
\end{equation}
where $s_j = \alpha + \imath \beta \cot \theta_j$ and $w_j = \frac{2\beta}{\sin^2 \theta_j}$, for $\theta_j \in (0,\pi)$. \footnote{A similar approximation to the minimization problem can be  obtained using the integration quadrature scheme in \cite{boyd_1987} in the infinite interval integral in $\norm{f}^p_{H^p(\IC_\alpha,\IV)}$.}


\subsection{Real-valued and complex-valued bases} \label{ssec:real-basis}
We present a set of theorems that establish additional properties of the proposed method, specifically analyticity and phase shift. More detailed explanation of the theorems used can be found in \cite{rudin_real_1987,king_hilbert_2009_1,titchmarsh_introduction_1986}.
\begin{definition}[Hilbert transform] \label{df:Hilbert}
    The Hilbert transform of $f$ can be thought of as the convolution of $f(t)$ with the function $h(t) = \frac{1}{\pi t}$, known as the Cauchy kernel. This represents a phase-shift of $\frac{\pi}{2}$ to the components of the function $f$. The Hilbert transform of a function  $f(t)$ is given by
    \begin{equation}\label{eq:htrans}
        \mathcal{H} f(x) = \frac{1}{\pi} \dashint_{-\infty}^{\infty} \frac{f(t-\tau)}{t} \odif{t} = \frac{1}{\pi} \lim_{\epsilon \rightarrow 0^+} \int_{|t|> \epsilon} \frac{f(t-\tau)}{t} \odif{t}
    \end{equation}
\end{definition}

\begin{theorem}[Titchmarsh] \label{th:Titchmarsh}
Let $\widehat{f} \in L^2(\imath\IR)$. If $\widehat{f}(s)$ satisfies any one of the following four conditions, then it satisfies all four conditions. The real and imaginary parts of $\widehat{f}(s)$, $\Re \{ \widehat{f}(s) \}$, and $\Im \{\widehat{f}(s)\}$, respectively, satisfy the following:
    \begin{enumerate} [label=(\roman*)]
    \itemsep0em 
    \item $\Im \widehat{f}(s) = \mathcal{H} \Re \widehat{f}(s)$;
    \item $\Re \widehat{f}(s) = -\mathcal{H} \Im \widehat{f}(s)$;
    \item If $f(t)$ denotes the inverse Fourier-Laplace transform \footnote{Note that the Laplace transform is equivalent to the Fourier transform with a change of variable for $f \in \IR_+$.} of $\widehat{f}(s)$, then $f(t)=0$, for $t<0$;
    \item $\widehat{f}(\sigma + \imath \omega)$ is an analytic function in the upper half plane and, for almost all $x$, $\widehat{f}(\sigma) = \lim_{y \rightarrow 0+} \widehat{f}(\sigma + \imath \omega)$.
    \end{enumerate} 
\end{theorem}

\begin{proposition}{} \label{prop:HL2}
    For $f \in L^2(\IR)$, $\widehat{f} = \mathcal{L} \{f\}$, the Hilbert transform $\mathcal{H}f$ is bounded in $L^2(\IR)$ as
    \begin{equation} \label{eq:norm-hilbert}
        \norm{\mathcal{H}f}_{L^2(\IR)} = \norm{\widehat{\mathcal{H}f}}_{L^2(\imath \IR)} = \norm{\widehat{f}}_{L^2(\imath \IR)} =\norm{f}_{L^2(\IR)}.
    \end{equation}
\end{proposition}

\begin{proof}
    We present here a informal proof similar to the one in \cite{taoLECTURENOTES247A}. A more formal proof using the Sochocki-Plemej theorem can be found in \cite{taoLECTURENOTES247A,chaudhury_lp_2012}. We write the Hilbert transform as the convolution $\mathcal{H}f = f * \frac{1}{\pi t}$ and by applying the Laplace transform and the convolution theorem we write
    \begin{equation}
        \widehat{\mathcal{H}f} (s) = \widehat{f}(s) \widehat{\frac{1}{\pi t}}(s) =-\imath \sgn(s)\widehat{f}(s).
    \end{equation}
    for $s \in \imath \IR$.
    Now, by using the Plancherel-Parseval identity we get
    \begin{equation}
        \int \abs*{\widehat{\mathcal{H}f}(s)}^2\odif{s} = \int \abs*{\widehat{f}(s)}^2\odif{s} = \int \abs*{f(t)}^2\odif{t}.
    \end{equation}
\end{proof}

\begin{corollary} \label{cor:H-rot}
    For $\widehat{f} \in H^2(\IC_\alpha)$ and $f=\mathcal{L}^{-1} \{\widehat{f} \} \in L^2(\IR_+)$, the norm of the limit function $\widehat{f}_\alpha$ can be defined using only its real part as
    \begin{align}
        \norm*{ \widehat{f}_{\alpha} }_{L^2(\imath\IR,\IV)} = \sqrt{2} \norm*{ \Re \{\widehat{f}_{\alpha} \} }_{L^2(\imath\IR,\IV)}
    \end{align}
    
\end{corollary}
\begin{proof}
    Noting that the function is only defined for $t>0$, we verify that the third condition in \cref{th:Titchmarsh} is fulfilled and therefore the other three. Also using \cref{prop:HL2} and then writing $\widehat{f}_\alpha$ as its real and imaginary parts, we get
    \begin{align}
        \norm*{\widehat{f}_{\alpha}}_{L^2(\imath\IR,\IV)} &= \int_{-\infty}^{\infty} \norm*{ \Re \{\widehat{f}_\alpha(s)\} + \imath \mathcal{H}\Re \{\widehat{f}_\alpha(s) \}}^2_{\IV} \odif{s} \\
        &= \int_{-\infty}^{\infty} \left( \norm*{ \Re \{\widehat{f}_\alpha(s)\}} + \norm*{\mathcal{H}\Re \{\widehat{f}_\alpha(s) \}}^2_{\IV} \right) \odif{s} = \int_{-\infty}^{\infty} \left( \norm*{ \Re \{\widehat{f}_\alpha(s)\}} + \norm*{\Re \{\widehat{f}_\alpha(s) \}}^2_{\IV} \right) \odif{s}
    \end{align}
\end{proof}

Using the results in \cref{cor:H-rot} we write a version of \cref{eq:fPOD} where we compute the basis using only the real part of the frequency solutions. 
\begin{equation} \label{eq:fPOD-real}
    \Phi =  \argmin_{\substack{\Phi \in \IR^{n_h \times n_r}: \\ \Phi^{\top} \Phi = \bm{I}_{n_r \times n_r}}}
    \sum_{j=1}^{n_s} w_j \norm*{\Re \{\widehat{\bm{\mathsf{u}}}(s_j) \} - \sum_{k=1}^{n_r} \dotp{\phi_k}{\Re \{ \widehat{\bm{\mathsf{u}}}(s_j) \}} \phi_k}^2_{L^2(\IT,\IV_h)}
\end{equation}

While analyticity ensures that the real and imaginary parts of \( \hat{u}(s) \) are functionally linked via the Hilbert transform, this relationship does not imply that the real parts alone span the same space as the full complex snapshots. Specifically, the set \( \{\Re[\hat{u}(s_j)]\} \) spans a real subspace of \( \IR^{n} \), which is strictly contained in the complex span \( \text{span}_\IC\{\hat{u}(s_j)\} \). As a result, \gls{pod} applied only to real parts may require more modes to reach the same approximation quality as a complex or symplectically lifted basis. In \cref{ssec:complexified-basis} we further explore this distinction and in \cref{th:symplectic} we present a symplectic lift method to construct the reduced-basis, which spans a real vector space isomorphic to the original complex subspace.

\subsection{Orthogonal, complex and symplectic structure of the basis} \label{ssec:complexified-basis}

In this part we want to further discuss the results from \cref{ssec:real-basis} using the scope of differential geometry.
Let us write the reduced solution in the frequency domain solution space $\widehat{u}_h^{\textrm{R}} \in \IV_h^{\IC,\textrm{R}}$ in terms of its real and imaginary parts as $\widehat{u}_h^{\textrm{R}} = \Re \{ \widehat{u}_h^{\textrm{R}} \} + \imath \Im \{ \widehat{u}_h^{\textrm{R}} \}$. Since in this section we work with both real and complex-valued spaces $\IV_h$, we indicate them as $\IV_h^{\IR}$ and $\IV_h^{\IC}$. Using the results from \cref{th:Titchmarsh,prop:HL2}, let us also consider the following matrix representation of the complex-valued reduced solution
\begin{equation} \label{eq:complex-matrixform}
    \widehat{u}_h \approx \Phi^\IC \widehat{u}_h^{\textrm{R}} =
    \begin{bmatrix}
        \Phi^\Re & -\Phi^\Im \\
        \Phi^\Im & \Phi^\Re
    \end{bmatrix}
    \begin{bmatrix}
        \Re \{ \widehat{u}_h \}  \\
        \Im \{ \widehat{u}_h \}
    \end{bmatrix} =
    \begin{bmatrix}
        \Phi^\Re & -\Phi^\Im \\
        \Phi^\Im & \Phi^\Re
    \end{bmatrix}
    \begin{bmatrix}
        \Re \{ \widehat{u}_h \}  \\
        \mathcal{H} \Re \{ \widehat{u}_h \}
    \end{bmatrix},
\end{equation}
with $\Phi^\Re$ and $\Phi^\Im$ the real and imaginary parts of the complex-valued basis $\Phi^\IC$. Note that by using \cref{eq:fPOD-circle} we construct the complex-valued basis $\Phi^\IC$, as we do it by solving an \gls{svd} of a complex-value data set. As we discussed in this section, and as a direct consequence of \cref{th:PW,cor:isometry} the reduced basis we are constructing is a basis for both the frequency-domain and the time-domain problems. Thus, we also write the real-value reduced solution for the time domain as
\begin{equation} \label{eq:real-matrixform}
    u_h \approx \Phi^\IC u_h^{\textrm{R}} =
    \begin{bmatrix}
        \Phi^\Re & -\Phi^\Im \\
        \Phi^\Im & \Phi^\Re
    \end{bmatrix}
    \begin{bmatrix}
        \Re \{ u_h^{\textrm{R}} \}  \\
        \Im \{ u_h^{\textrm{R}} \}
    \end{bmatrix} =
    \begin{bmatrix}
        \Phi^\Re & -\Phi^\Im \\
        \Phi^\Im & \Phi^\Re
    \end{bmatrix}
    \begin{bmatrix}
        \Re \{ u_h^{\textrm{R}} \}  \\
        \mathcal{H} \Re \{ u_h^{\textrm{R}} \}
    \end{bmatrix}.
\end{equation}
Since we know that $u_h^{\textrm{R}}$ in the time-dependent problem is real-valued, we get
\begin{equation} \label{eq:realfromcomplex}
    \Phi^\Im \Re \{ u_h^{\textrm{R}} \} + \Phi^\Re \Im \{ u_h^{\textrm{R}} \} = \Phi^\Im \Re \{ u_h^{\textrm{R}} \} + \Phi^\Re \mathcal{H} \Re \{ u_h^{\textrm{R}} \} = 0.
\end{equation}

This matrix representation implies that the complex $n_s$-dimensional space $\IV^\IC_h$, can be defined using a real space $\IV^\IR_h$. To explain this let us start by writing some basic definitions and theorems used in holomorphic differential geometry, for which we follow \cite[chapter~4]{lewisNotesGlobalAnalysis} and \cite{maninLinearAlgebraGeometry1989}.
\begin{definition} [$\mathcal{H}$ as the linear complex structure of $\IV_h^\IR$]
    We define a linear complex structure $\mathcal{J}$ on a real vector space $\IV^\IR$ as the endomorphism $\mathcal{J}: \IV_h^\IR \to \IV_h^\IR$ such that $\mathcal{J} \cdot \mathcal{J} =- I$, with $I$ the identity.
    Noting that by definition the Hilbert transform is an anti-involution $\mathcal{H}^2 = -I$, we define the endomorphism
\begin{equation} \label{eq:J-complex}
    \mathcal{J}_{\mathcal{H}} = 
    \begin{bmatrix}
        0 & -\mathcal{H} \\
        \mathcal{H} & 0
    \end{bmatrix}
\end{equation}
as a linear complex structure on $\IV_h^\IR$.
\end{definition}

\begin{definition} [Real structure of $\IV_h^\IC$] \label{def:realstruc}
    A real structure on a complex vector space $\IV^\IC$ is an antilinear map $\varsigma: \IV_h^\IC \to \IV_h^\IC$ such that $\varsigma \cdot \varsigma = I$. Let us write any $u \in \IV_h^\IC$ as the sum of $\widehat{u} = \widehat{u}_\Re + \widehat{u}_\Im = \frac{1}{2} (\widehat{u} + \varsigma \widehat{u}) + \frac{1}{2} (\widehat{u} - \varsigma \widehat{u})$. This way, we define the subspaces
    \begin{align}
        \IV_h^\Re = \{\widehat{u} \in \IV_h^\IC | \varsigma \widehat{u} = \widehat{u}^* \}, \quad \text{and} \quad
        \IV_h^\Im = \{\widehat{u} \in \IV_h^\IC | \varsigma \widehat{u} = -\widehat{u}^* \}
    \end{align}
    with $\widehat{u}^*$ the complex conjugate of $\widehat{u}$. Note that since $\widehat{u}_\Re$ and $\widehat{u}_\Im$ are orthogonal, then $\dim_\IR \IV_h^\Re = \dim_\IR \IV_h^\Im = \dim_\IC \IV_h^\IC$.
\end{definition}

\begin{definition} [Complexification and realification] \label{def:complexification}
    For a real vector space $\IV_h^\IR$ we may define its complexification by extension of scalars and for a complex vector space $\IV_h^\IC$ we may define its realification by the restriction of scalars. We get the relation $\IV_h^{\IC} = \IC \otimes_{\IR} \IV_h^\IR$.
\end{definition}

\begin{proposition} [Basis representation of $\IV_h^\IC$ and $\IV_h^\IR$] \label{prop:basisrc}
    Let $\{\varphi_1,\ldots,\varphi_{n_s}\}$ be a basis of $\IV_h^\IC$. Then, $$\{\varphi_1,\ldots,\varphi_{n_s}, \mathcal{J}_\mathcal{H}(\varphi_1),\ldots,\mathcal{J}_\mathcal{H}(\varphi_{n_s})\}$$ forms a real basis for $\IV_h^\IR$. Note that $\dim_\IR \IV_h^\IR = 2\dim_\IC \IV_h^\IC$.
\end{proposition}

\begin{proposition}[Totally real subspace] \label{prop:realsubs}
    Let $\IV_h$ be a finite-dimensional $\IR$-vector space with linear complex structure $\mathcal{J}_{\mathcal{H}}$. A subspace $\IV_h^{\IR} \subseteq \IV_h$ is totally real if $\mathcal{J}_{\mathcal{H}}(\IV_h^{\IR}) \cap \IV_h^{\IR} = \{0\}$. The next two conditions are equivalent 
    \begin{itemize}
        \item $\IV_h^{\IR}$ is a totally real subspace.
        \item $\mathcal{J}_{\mathcal{H}}(\IV_h^{\IR})$ and $\IV_h^{\IR}$ are orthogonal in the inner product sense.
    \end{itemize}
\end{proposition}

\begin{corollary} \label{cor:j-basis}
    Let $\IV_h^{\textrm{R}}$ be a reduced space for the time dependent problem and let $\{\phi_1,\dots,\phi_{n_s}\}$ be a real-valued reduced basis of $\IV_h^\IR$, then $\{\phi_1,\dots,\phi_{n_s},\mathcal{J}_{\mathcal{H}}(\phi_i),\ldots,\mathcal{J}_{\mathcal{H}}(\phi_{n_r})\}$ is a basis of $\IV_h^\IC$.
\end{corollary}
\begin{proof}
    Let $v^{\textrm{R}} \in L^2(\IT,\IV_h^\IR)$. Then, to prove the orthogonality between $\phi_j$ and $\mathcal{J}_{\mathcal{H}}(\phi_j)$ we can show that any $v^{\textrm{R}}_j$ is orthogonal to $\mathcal{J}_{\mathcal{H}}(v^{\textrm{R}}_j)$. Noting that $\mathcal{H}$ is anti-self adjoint we write
    \begin{equation}
        \dual{\mathcal{H}(v_j)}{v_j} = \dual{\mathcal{H}^2(v_j)}{\mathcal{H}(v_j)} = \dual{v_j}{-\mathcal{H}(v_j)} = \dual{\mathcal{H}(v_j)}{v_j} =0.
    \end{equation}
    Using \cref{prop:realsubs} we see that $\IV_h^\IR$ is a totally real subspace and can be described with $\Phi^\Re$. We check that for \cref{eq:realfromcomplex} we satisfy $\Phi^\Im \Re \{ u_h^{\textrm{R}} \} + \Phi^\Re \Im \{ u_h^{\textrm{R}} \} =0$ is by assuming that $\Phi^\Re$ and $\Phi^\Im$ are orthogonal, and that $\Re \{ u_h^{\textrm{R}} \} \in \spn \{ \phi_1^\Re, \ldots, \phi_{n_r}^\Re \}$ and $\Im \{ u_h^{\textrm{R}} \} \in \spn \{ \phi_1^\Im, \ldots, \phi_{n_r}^\Im \}$.
\end{proof}

So far we have described real and complex vector spaces. Here it is clear that $\IV_h^\IC$ represents the frequency-domain solution space, but there is not an unique way to construct a real space or obtain a real-valued reduced solution. Moreover, we see that the spaces $\IV_h^\Re$ and $\IV_h^\IR$ are not equivalent.

This non-equivalence follows directly from the structure of the snapshot spaces. Given a set of complex-valued snapshots $\hat{u}(s_j) \in \mathbb{C}^n$, we define the complex reduced space $V_h^{\mathbb{C}} := \operatorname{span}_{\mathbb{C}}\{\hat{u}(s_j)\}$, and its realification corresponds to a real subspace of $\mathbb{R}^{2n}$. In contrast, the space $V_h^{\Re} := \operatorname{span}_{\mathbb{R}}\{\Re\{\hat{u}(s_j)\}\}$ lies in $\mathbb{R}^n$ and does not include the imaginary components. Although the imaginary parts $\Im\{\hat{u}(s_j)\}$ are functionally related to the real parts via the Hilbert transform due to analyticity (as noted in \cref{cor:H-rot}), this relation is non-linear and does not preserve span equivalence in the finite-dimensional \gls{pod} setting.

Consequently, the reduced space $V_h^{\Re}$ cannot reproduce the full solution space captured by $V_h^{\mathbb{C}}$, even though it may yield an effective approximation in practice. To recover a real-valued basis that faithfully represents the original complex space, we introduce the symplectic lifting approach in the next section, which constructs a real space in $\mathbb{R}^{2n}$ that is isomorphic to the complex snapshot span.

\begin{definition}[Symplectic lift] \label{def:symlift}
Given a complex vector $\hat{u} \in \mathbb{C}^n$, we define its \emph{symplectic lift} as the real vector in $\mathbb{R}^{2n}$ given by
\begin{equation}
    \mathcal{L}_S(\hat{u}) := \begin{bmatrix} \Re\{\hat{u}\} \\ \Im\{\hat{u}\} \end{bmatrix}.
\end{equation}
We extend this definition component-wise to a set of snapshots $\{\hat{u}_j\}_{j=1}^{n_s}$ to obtain the real-valued snapshot matrix
\begin{equation}
    S_S := \left[ \mathcal{L}_S(\hat{u}_1), \ldots, \mathcal{L}_S(\hat{u}_{n_s}) \right] \in \mathbb{R}^{2n \times n_s}.
\end{equation}
\end{definition}

Using \cref{def:symlift,prop:basisrc} in $\widehat{\bm{\mathsf{S}}}$ (\cref{eq:fsnap}), and replacing $\mathcal{J}_{\mathcal{H}} (\Re \{ \widehat{\bm{\mathsf{u}}}_j \}) = \Im \{ \widehat{\bm{\mathsf{u}}}_j \}$ as per its definition. We obtain the rearranged snapshot set
\begin{equation} \label{eq:ssnap}
    \bm{\mathsf{S}}^S = \{ \Re \{ \widehat{\bm{\mathsf{u}}}_1 \}, \ldots, \Re \{ \widehat{\bm{\mathsf{u}}}_{n_s} \}, \Im \{ \widehat{\bm{\mathsf{u}}}_1 \}, \ldots, \Im \{ \widehat{\bm{\mathsf{u}}}_{n_s} \} \}
\end{equation}
Now, if we compute the \gls{svd} and obtain a basis for this rearranged set, we are effectively constructing a symplectic basis. This technique is the same as the one presented in \cite{peng_symplectic_2015} for Hamiltonian problems as the cotangent lift method.
\begin{theorem} [Cotangent lift method] \label{th:symplectic}
    Let $\widehat{\bm{\mathsf{u}}}$ be the frequency discrete solutions and $\bm{\mathsf{S}}^S$ the data set as described in \cref{eq:ssnap}. The solution of the minimization problem in \cref{eq:fPOD-circle} is given by $\bm{\mathsf{S}}^S = \Phi^S \Sigma (\Psi^S)^*$, which results in the symplectic basis
    \begin{equation}
        \Phi^S =
    \begin{bmatrix}
        \Phi & 0 \\
        0 & \Phi
    \end{bmatrix}
    \end{equation}
    which is the optimal projection of the sample set $\widehat{\bm{\mathsf{S}}}$ onto the column space of the symplectic vector space $\IR^{2n_s}$.
\end{theorem}
\begin{proof}
    Writing  \cref{eq:fPOD-circle} in terms of the real and imaginary parts, we get
    \begin{align}
        &\norm*{
        \begin{bmatrix}
            \Re \{ \widehat{\bm{\mathsf{u}}} \} \\
            \Im \{ \widehat{\bm{\mathsf{u}}} \}
        \end{bmatrix} -
        \begin{bmatrix}
            \Phi & 0 \\
            0 & \Phi
        \end{bmatrix}
        \begin{bmatrix}
            \Phi^\top & 0 \\
            0 & \Phi^\top
        \end{bmatrix}
        \begin{bmatrix}
            \Re \{ \widehat{\bm{\mathsf{u}}} \} \\
            \Im \{ \widehat{\bm{\mathsf{u}}} \}
        \end{bmatrix}
        }_{L^2} =
        \norm*{
        \begin{bmatrix}
            \left( I_{n_s} -\Phi \Phi^\top \right) \Re \{ \widehat{\bm{\mathsf{u}}} \} \\
            \left( I_{n_s} -\Phi \Phi^\top \right) \Im \{ \widehat{\bm{\mathsf{u}}} \}
        \end{bmatrix}
        }_{L^2} \\ \vspace{0.2cm}
        &\quad = \norm*{
        \begin{bmatrix}
            I_{n_s} -\Phi \Phi^\top
        \end{bmatrix}
        \begin{bmatrix}
            \Re \{ \widehat{\bm{\mathsf{u}}} \} &
            \Im \{ \widehat{\bm{\mathsf{u}}} \}
        \end{bmatrix}
        }_{L^2}
        = \norm*{
        \begin{bmatrix}
            \Re \{ \widehat{\bm{\mathsf{u}}} \} &
            \Im \{ \widehat{\bm{\mathsf{u}}} \}
        \end{bmatrix} -
        \Phi \Phi^\top
        \begin{bmatrix}
            \Re \{ \widehat{\bm{\mathsf{u}}} \} &
            \Im \{ \widehat{\bm{\mathsf{u}}} \}
        \end{bmatrix}
        }_{L^2}
    \end{align}
    giving the minimization
    \begin{equation} \label{eq:fPOD-sym}
        \Phi^S =  \argmin_{\substack{\Phi^S \in \IR^{n_h \times n_r}: \\ (\Phi^S)^{\top} \Phi^S = \bm{I}_{n_r \times n_r}}}
        \sum_{j=1}^{n_s} w_j \norm*{ \left[ \Re \{\widehat{\bm{\mathsf{u}}}(s_j) \}, \Im \{\widehat{\bm{\mathsf{u}}}(s_j) \} \right] - \sum_{k=1}^{n_r} \dotp{\phi_k^S}{\left[ \Re \{\widehat{\bm{\mathsf{u}}}(s_j) \}, \Im \{\widehat{\bm{\mathsf{u}}}(s_j) \} \right]} \phi_k^S}^2_{L^2(\IT,\IV_h)}
\end{equation}
\end{proof}

\begin{remark}
    The equivalence between the complex-valued basis, the symplectic basis, and the real-only \gls{pod} basis can be understood geometrically. The symplectic lift maps the complex vector space $\mathbb{C}^n$ into the real vector space $\mathbb{R}^{2n}$, endowing it with both a canonical symplectic form and a compatible complex structure. When the lifted snapshot matrix is used to compute a real basis, the resulting space naturally preserves both phase and amplitude information, just as the complex basis does.
    This triplet of structures ---orthogonal inner product, symplectic form, and complex structure--- defines a Kähler manifold in the reduced space. In this setting, the symplectic basis can be viewed as a real coordinate representation of the complex basis.
\end{remark}

Using the results and definitions from this section, we define different approaches to the construction of the reduced basis for the time-dependent problem.
\begin{itemize}
    \item (\Cref{th:PW,cor:isometry}). The complex-valued basis $\Phi_{\IC}$ obtained from \cref{eq:fPOD-circle}. This is equivalent to constructing a reduced space for $\IV_h^\IC$, which after \cref{cor:j-basis} we know results in a real-valued solution $u_h^{\textrm{R}}$
    \item (\Cref{th:Titchmarsh,cor:H-rot}). The real-valued basis $\Phi_\IR$ obtained from \cref{eq:fPOD-real}. This is equivalent as the reduced basis for a space $\IV_h^\IR$. Note that due to the smaller dimension of this space (\cref{prop:basisrc}) we can expect a worse performing reduced basis when using the same number of samples $n_s$.
    \item (\Cref{prop:realsubs,cor:j-basis}). The real-valued basis $\Phi_{\Re \IC}= \Re \{\Phi_{\IC} \}$ obtained by solving \cref{eq:fPOD-circle} and taking the real part of the resulting basis. This is equivalent as the reduced basis of the space $\IV_h^\Re$. This basis naturally solves a real-valued problem and spans a space of the same size as $\Phi_\IC$.
    \item (\Cref{th:symplectic}). The symplectic basis $\Phi_{S}$ obtained from \cref{eq:fPOD-sym}. We can easily obtain this basis by computing an \gls{svd} of the rearranged data set in \cref{eq:ssnap}. With this construction we get a real-valued reduced basis that spans a space of the same size as $\Phi_\IC$.
\end{itemize}

In \cref{ssec:Heat}, we show that all the different reduced bases adequately solve the time-dependent problem. We also observe that the bases $\Phi_{\mathbb{C}}$ and $\Phi_{S}$ exhibit similar behaviour, consistent with the equivalence relations discussed in this section.
\section{Numerical examples} \label{sec:Examples}
In this section, we present three numerical examples to evaluate the proposed method. The first example employs the wave equation to assess the method's convergence properties, computational efficiency, and model parameter sensitivity ($\beta$, $n_r$, $n_s$). The second example uses a heat equation to compare the different basis construction approaches explained in \cref{ssec:complexified-basis} and demonstrate their relative performances. In the third example, we apply the method to a convection-dominated problem to highlight the stabilizing properties of the frequency-domain approach when handling problems that typically exhibit numerical instabilities.

To numerically evaluate the reduced order solutions in all the problems, we define the $H^1(\setI \times \Omega)$-error as
\begin{equation} \label{eq:H1-error}
    H^1-\textrm{error}(u_h^{\textrm{R}}) = \norm{u_h - u_h^{\textrm{R}}}_{H^1(I \times \Omega)}^2 
    = \frac{1}{n_t} \sum_{j=1}^{n_t} \left( \norm{u_h(t_j) - u_h^{\textrm{R}}(t_j)}_{H^1(\Omega)} \right)^2,
\end{equation}
and the relative $H^1(\setI \times \Omega)$-error as
\begin{equation} \label{eq:H1-error-rel}
    \textrm{relative } H^1-\textrm{error}(u_h^{\textrm{R}}) = 
    \frac{\norm{u_h - u_h^{\textrm{R}}}_{H^1(I \times \Omega)}}{\norm{u_h}_{H^1(I \times \Omega)}}.
\end{equation},
where $u_h$ represents the finite element solution and $u_h^{\textrm{R}}$ represents the reduced-order solution.

\subsection{Wave equation: Model parameters} \label{ssec:Wave}

\begin{figure}[H]
	\centering
    \begin{tblr}
    {cell{1}{3} = {r=2}{},}
	\begin{subfigure}[htb!]{0.3685\textwidth}
	    \centering
	    \includegraphics[width=\textwidth]{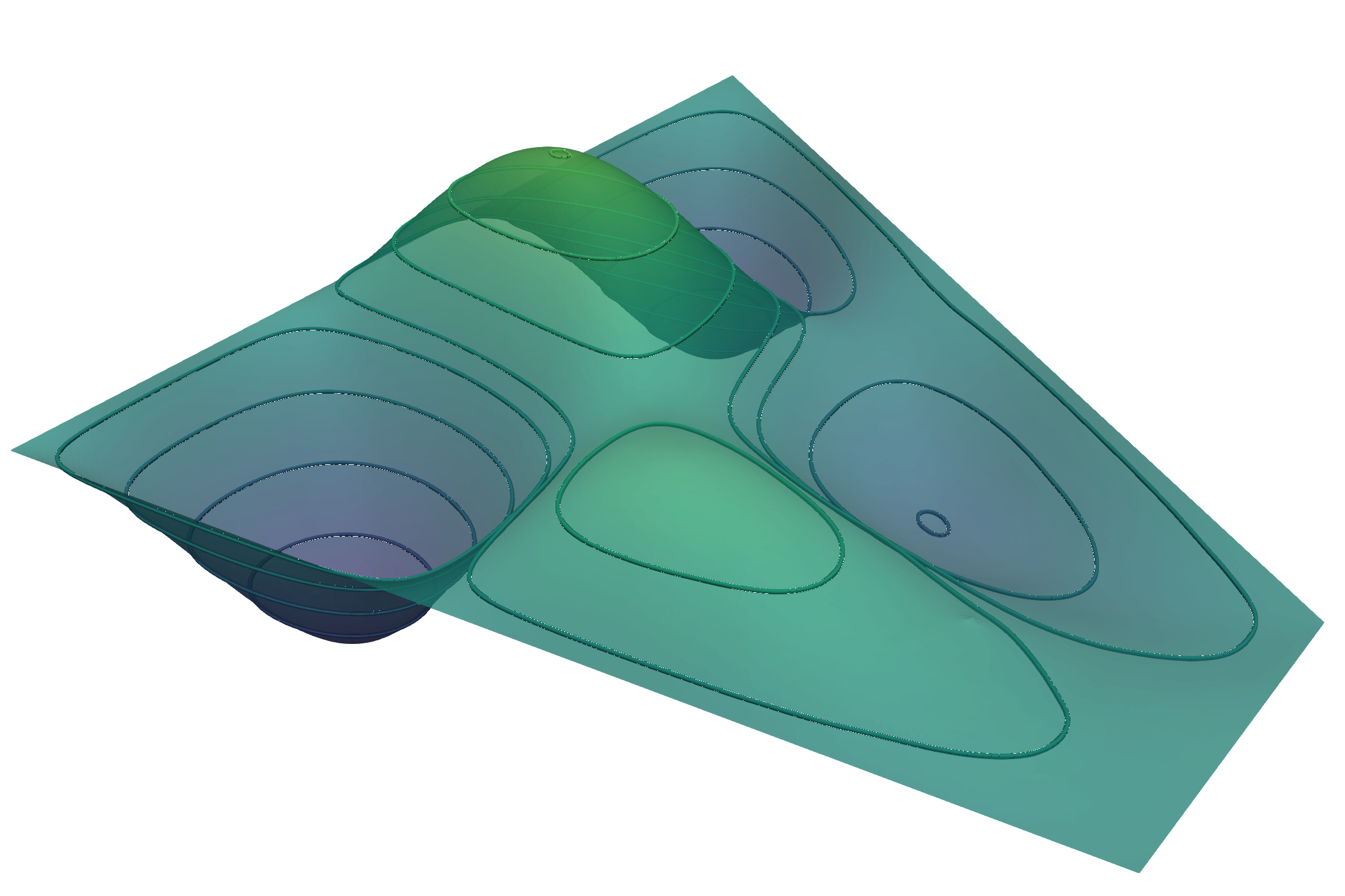}
        \caption{$t=23.5$} 
    \end{subfigure} &
    \begin{subfigure}[htb!]{0.3685\textwidth}
	    \centering
	    \includegraphics[width=\textwidth]{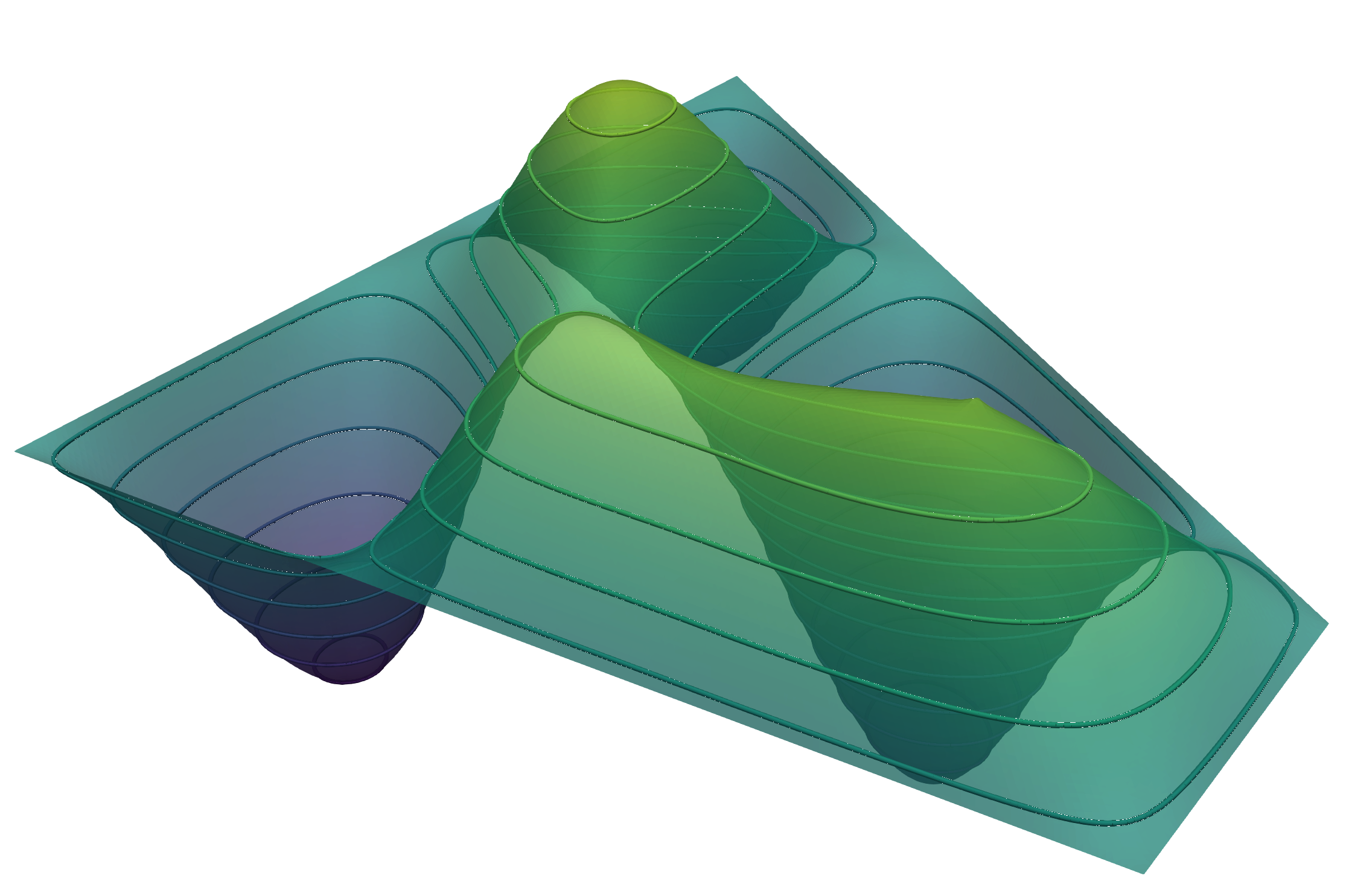}
        \caption{$t=24.0$}
    \end{subfigure} & \includegraphics[scale=0.16]{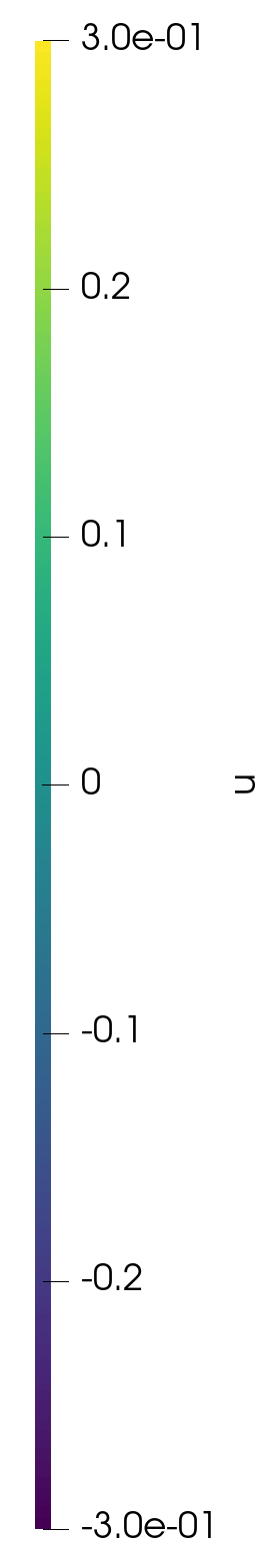} \\
    \begin{subfigure}[htb!]{0.3685\textwidth}
	    \centering
	    \includegraphics[width=\textwidth]{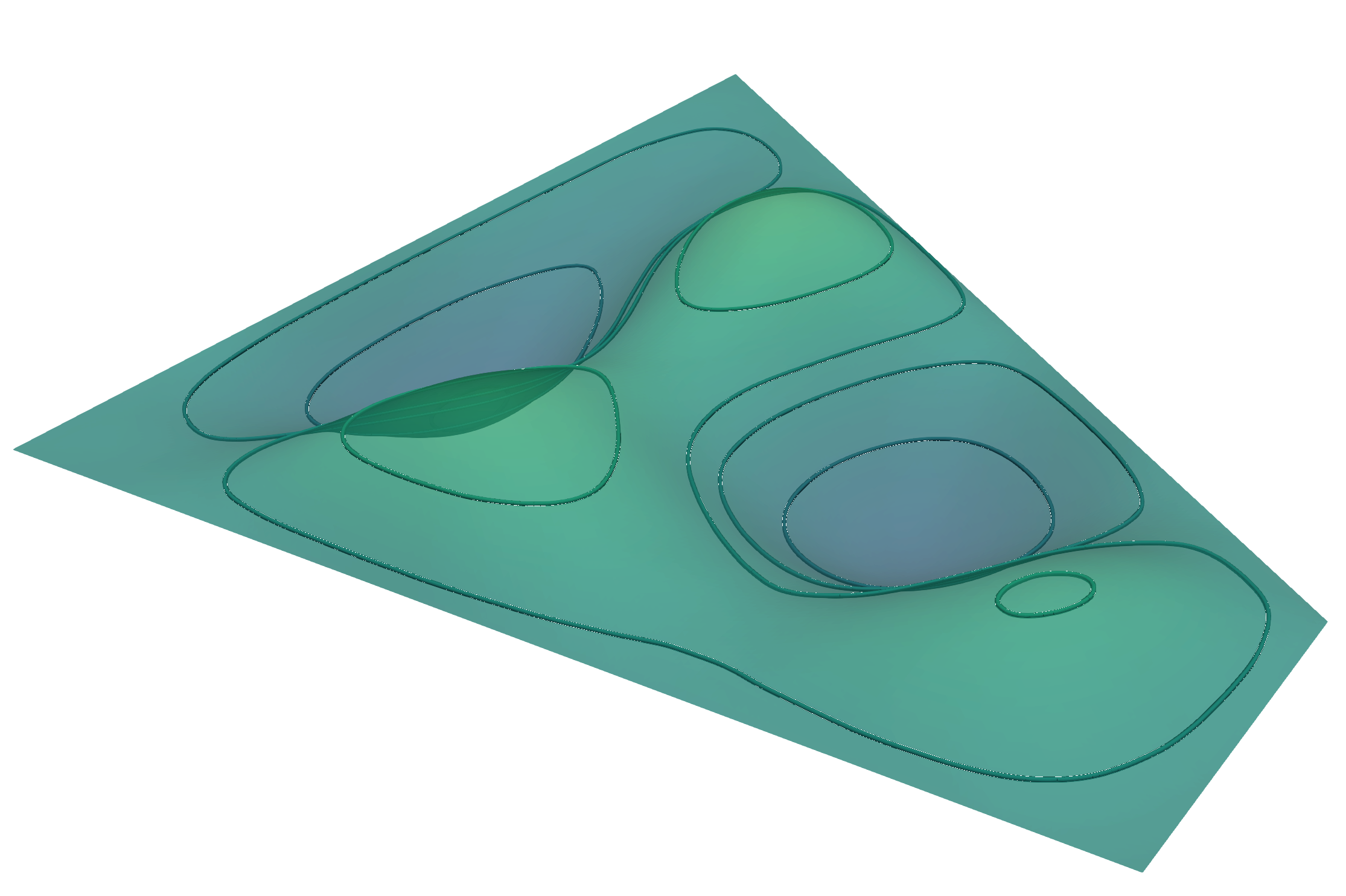}
        \caption{$t=24.5$} 
    \end{subfigure} &
    \begin{subfigure}[htb!]{0.3685\textwidth}
	    \centering
	    \includegraphics[width=\textwidth]{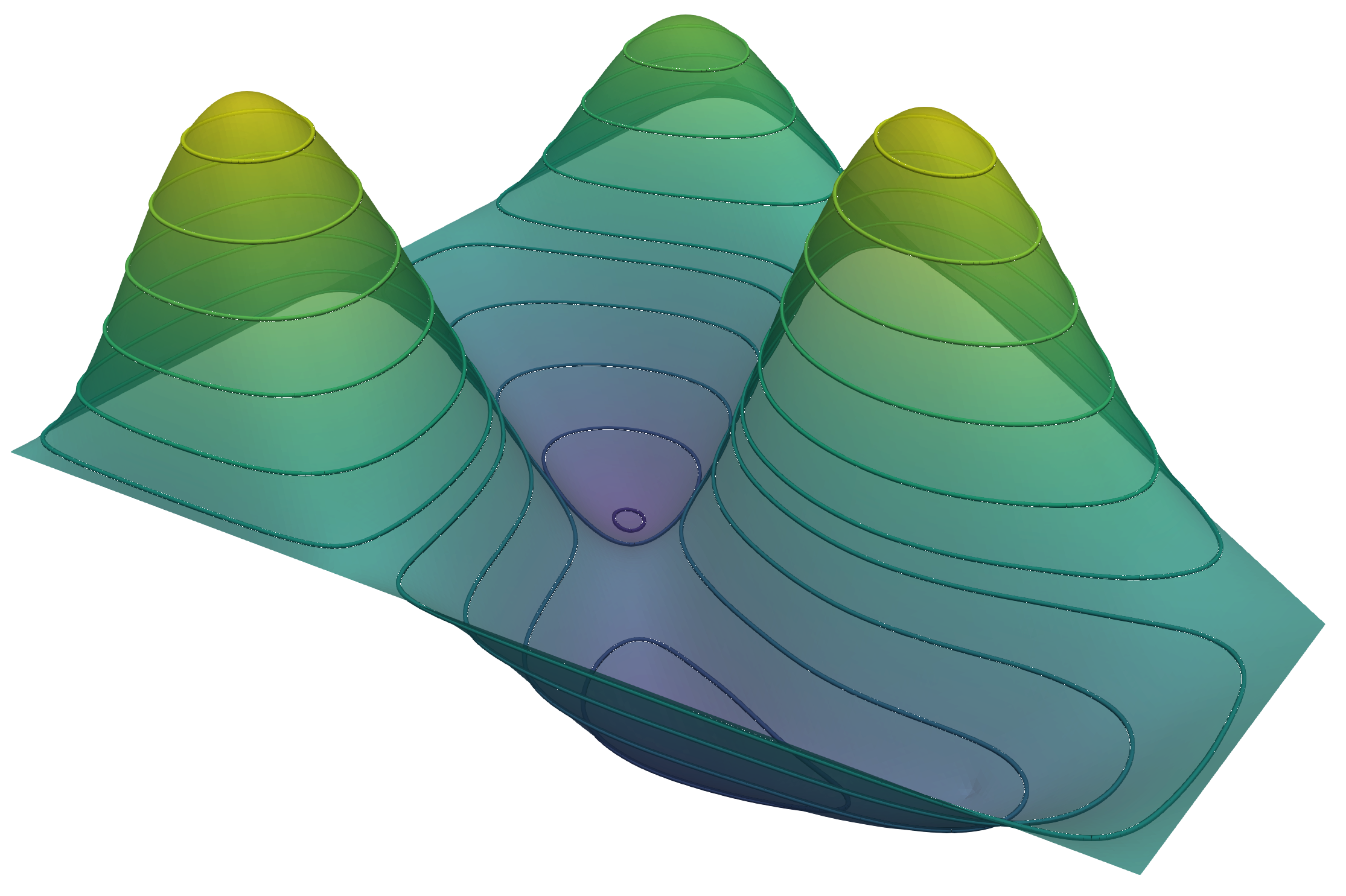}
        \caption{$t=25.0$}
    \end{subfigure} &
    \end{tblr}
	\caption{Wave problem: solution contours at selected time steps.} 
	\label{fig:wave-contour}
\end{figure}

As a first numerical example we solve the non-homogeneous wave equation
\begin{alignat}{2} \label{eq:wave}
    u_{tt} - c^2 \nabla^2 u &= b &\qquad \text{on} \; \Omega \times \setI_h \\
    u(\bm{x},t) &= u_D(\bm{x},t) = 0 &\qquad \text{on} \; \partial \Omega \times \setI_h \\
    u(\bm{x},0) &= u_0(\bm{x}) = 0 &\qquad \text{on} \; \Omega.
\end{alignat}
The domain $\Omega$ consists in a polygonal shape described by the vertices $[(0.0,0.0),(1.0,0.0),(1.0,0.3),(0.0,1.0)]$, over the time interval $\setI = (0,25]$. We discretize the problem using a mesh of $61204$ triangular elements and a time step $\delta t = 0.005$. We define the load with its components $b_x$ and $b_t$, and the load in the frequency domain $\widehat{b}$ as
\begin{align}
    b_x(\bm{x}) &= \exp \left(-\gamma^2 \left[(x_0-0.8)^2 + (x_1-0.15)^2 \right] \right), \\
    b_t(t) &= A \sin(\omega t) \e^{-\tau t}, \\
    \widehat{b}(s) &= \frac{\omega}{\omega^2 +(s+\tau)^2},
\end{align}
with $\bm{x}=[x_0,x_1]$, $\gamma=500$, $\omega=10$, $\tau=0.1$, $c=0.75$ and $A=10^4$.

The aim of this numerical example is to evaluate the behaviour of the solution for different values of $\beta$, and for different sizes of the sampling set $n_s$ and the reduced basis $n_r$. We follow the symplectic basis construction $\Phi_{S}$, fixing $\alpha = 1$. In \cref{fig:wave-contour} we show a contour plot of the solution at different time steps, here we can see the complexity of the solution and how unpredictable is its behaviour caused by the interference of the waves in the enclosed domain.

\begin{figure}[H]
	\centering
	\begin{subfigure}[htb!]{0.49\textwidth}
	    \centering
        \includegraphics{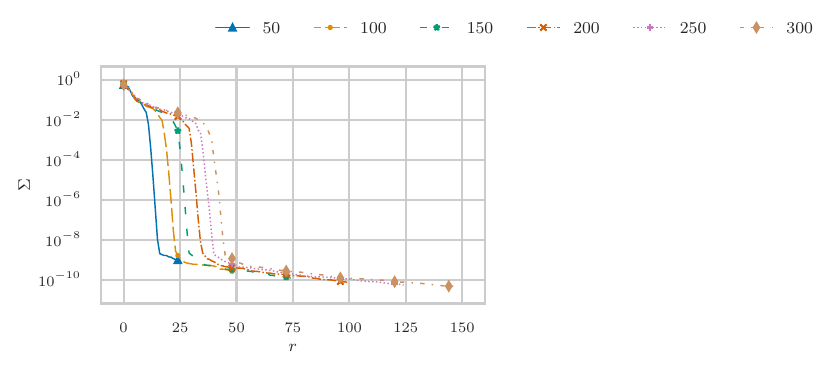}
        \caption{$\beta=2$.} 
        \label{fig:wave-basis-b1} 
    \end{subfigure}
    \begin{subfigure}[htb!]{0.49\textwidth}
	    \centering
	    \includegraphics{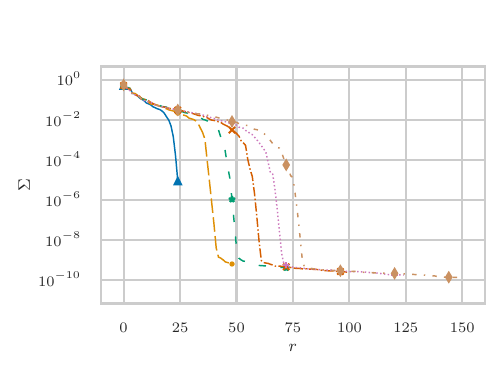}
        \caption{$\beta=5$.} 
        \label{fig:wave-basis-b2} 
    \end{subfigure}
    \begin{subfigure}[htb!]{0.49\textwidth}
	    \centering
	    \includegraphics{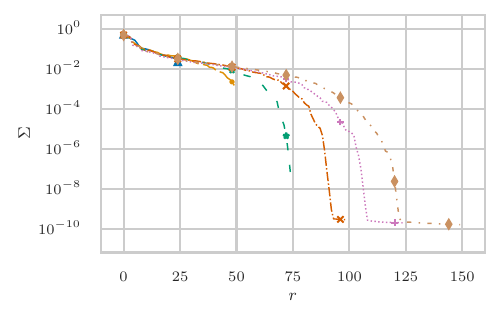}
        \caption{$\beta=10$.} 
        \label{fig:wave-basis-b3} 
    \end{subfigure}
    \begin{subfigure}[htb!]{0.49\textwidth}
	    \centering
	    \includegraphics{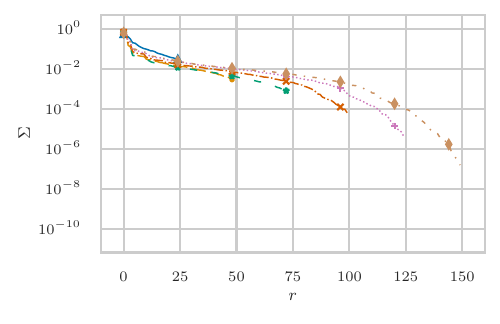}
        \caption{$\beta=15$.} 
        \label{fig:wave-basis-b4} 
    \end{subfigure}
    \caption{Wave problem: singular value decay for different values of \( \beta \).} 
	\label{fig:wave-basisS} 
\end{figure}

In \cref{fig:wave-basisS} we show the singular value decay $\Sigma$ for sample sizes $n_s = [50, 100, 150, 200, 250, 300]$ and $\beta = [2, 5, 10, 15]$. As expected, the singular value decay rate depends on the sampling set size $ns$, with larger sample sets providing more information about the solution behaviour. This decay also changes significantly with variations in the parameter $\beta$, likely because increasing $\beta$ shifts the sampling farther from $0$ along the line $\imath \IR$. As a result, the sampling set includes higher frequencies, which require more basis vectors for accurate representation.

\begin{figure}[H]
	\centering
    \begin{subfigure}[htb!]{0.49\textwidth}
     	\centering
        \includegraphics{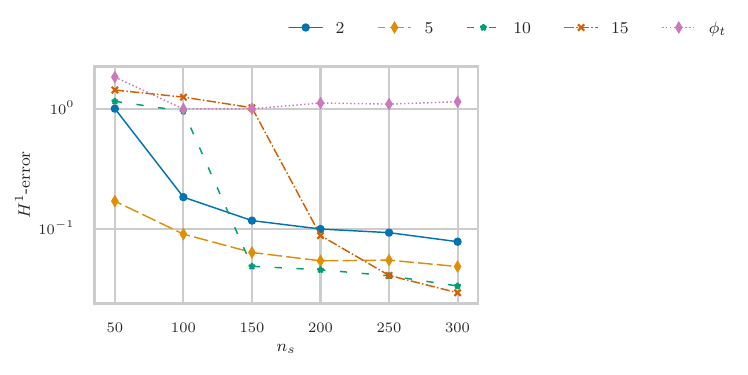}
        \caption{Error convergence for $n_s$ with $n_r=\frac{n_s}{2}$ .} 
	    \label{fig:wave-error-ns}
    \end{subfigure}
    \begin{subfigure}[htb!]{0.49\textwidth}
     	\centering
        \includegraphics{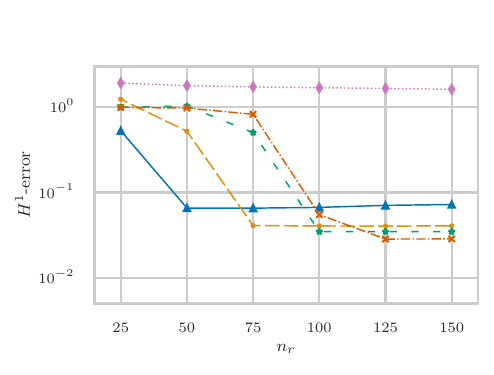}
        \caption{Error convergence for $n_r$ with $n_s=200$ samples.} 
	    \label{fig:wave-error-nr}
    \end{subfigure}
    \caption{Wave problem: relative \( H^1 \)-error for different values of \( \beta \).}
	\label{fig:wave-error}
\end{figure}

In \cref{fig:wave-error} we compare the relative error (\cref{eq:H1-error-rel}) for different values of $\beta$. \Cref{fig:wave-error-ns} displays the error convergence for sample size $n_s$, while keeping the basis size $n_r=\frac{n_s}{2}$. \Cref{fig:wave-error-nr} shows the error convergence for the reduced basis size $n_r$, fixing the sample size to $n_s=200$. In both figures we also compare a solution obtained with the standard \gls{pod}, by sampling the first $n_s=200$ finite element solutions of the time interval.

Although error decreases for bigger basis size in all the cases, we find different behaviour for basis computed with different constant $\beta$. As we increase the value of $\beta$ we are sampling higher frequencies and consequently we need a bigger amount of samples to appropriately capture the behaviour of the problem. This can be seen by comparing \cref{fig:wave-basisS,fig:wave-error}, where the point where the singular values begin to decay matches the decrease of the error. On the other hand, when we have enough samples as in \cref{fig:wave-error-ns}, the error decay slope is similar for all the cases, still requiring a bigger basis size for a bigger $\beta$. We should also note that the behaviour of the beta is not linear, as the case with $\beta=5$ has better approximation that the others.

\begin{figure}[H]
	\centering
    \includegraphics{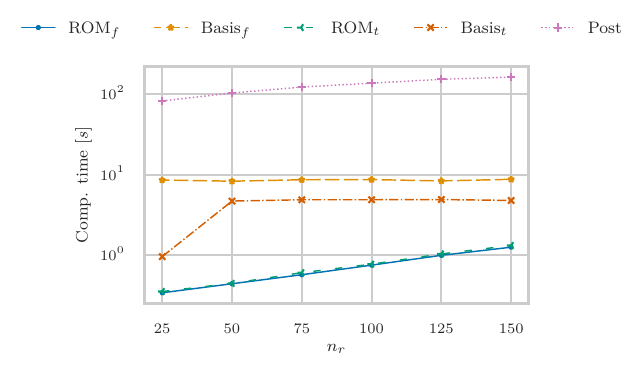}
    \caption{Wave problem: computation time for reduced-order solution and basis construction.}
	\label{fig:wave-time}
\end{figure}

We compare the computational cost of the \gls{frb} method and the traditional \gls{pod} method using $\beta=5$ and the same sample size $n_s=1000$. In \cref{fig:wave-time} we show the computational time of the construction of the basis, the reduced solution, and the post-process for different sizes of the basis. The finite element problem is solved in $1248$ seconds, and the sampling is perform in $212$ seconds. In this case we solve the frequency problem in a parallel way using 24 cores.

\subsection{Heat equation: Topology} \label{ssec:Heat}

\begin{figure}[H]
    \centering
	\begin{subfigure}[htb!]{0.49\textwidth}
	    \centering
	    \includegraphics{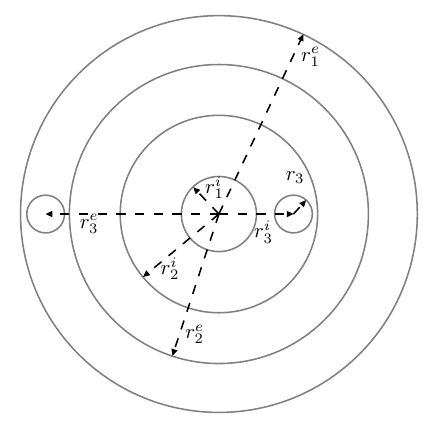}
	    \caption{Geometrical parameters.}
        \label{fig:dif-problem-1}
    \end{subfigure}
    \begin{subfigure}[htb!]{0.49\textwidth}
	    \centering
	    \includegraphics[width=0.9\textwidth]{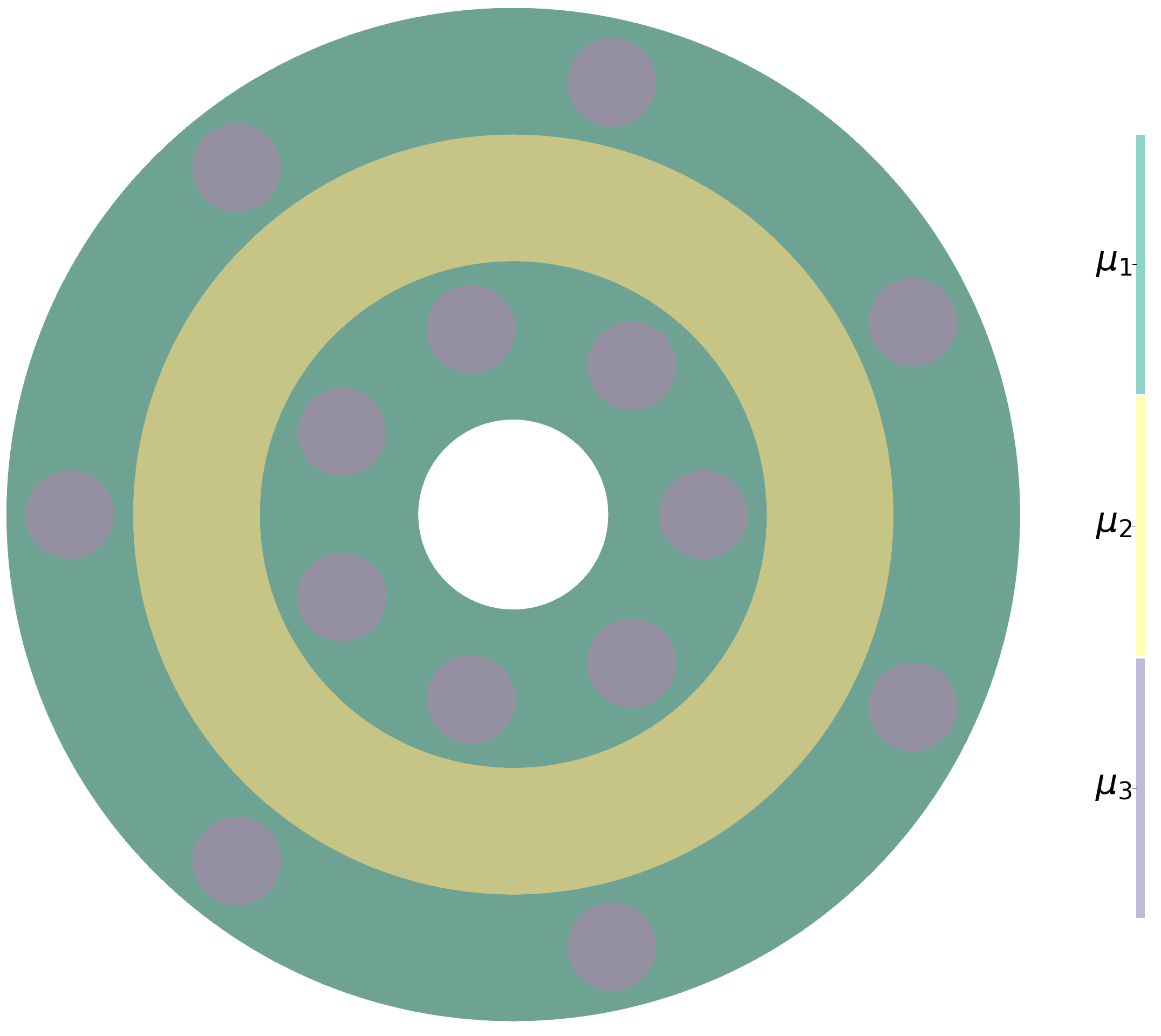}
	    \caption{Material properties ($\mu$).}
        \label{fig:dif-problem-2}
    \end{subfigure}
    \caption{Diffusion problem: domain geometry and material distribution.}
    \label{fig:dif-problem}
\end{figure}

As the second numerical example we solve a heat equation in a 2D circular domain over a time interval of $\setI=(0,100]$. The domain $\Omega$ consist of three successive rings, with seven circular inclusions in the inner and outer rings. The inclusions are distributed evenly in each ring shifted by $\pi$ with respect to each other. \Cref{fig:dif-problem} shows a schematic representation of the geometry and the physical properties of the problem. \Cref{tab:dif-data} shows the geometrical and physical parameters of the problem. We discretize the problem using a mesh of $55359$ triangular elements and a time step $\delta t = 0.1$.

\begin{table}[!htb]
\centering
\caption{Geometrical parameters and material properties of the diffusion problem.}
\begin{subtable}[ht]{0.4\textwidth}
    \centering
    {\small
    \begin{tabular}{ccccccc}
        \toprule
        $r_3$ & $r_3^e$ & $r_3^i$ & $r_1^e$ & $r_1^i$ & $r_2^e$ & $r_2^i$  \\ \midrule  
        0.35 & 3.5 & 1.5 & 4.0 & 0.75 & 3.0 & 2.0 \\ \bottomrule
    \end{tabular}}
    \label{tab:dif-geometry}
\end{subtable}
\begin{subtable}[ht]{0.3\textwidth}
    \centering
    {\small
    \begin{tabular}{ccc}
        \toprule
        $\mu_1$ & $\mu_2$ & $\mu_3$  \\ \midrule  
        0.075 & 0.5 & 27.2 \\ \bottomrule
    \end{tabular}}
    \label{tab:dif-material}
\end{subtable}
\label{tab:dif-data}
\end{table}

\begin{figure}[H]
	\centering
	\begin{subfigure}[htb!]{0.5\textwidth}
        \centering
        \includegraphics{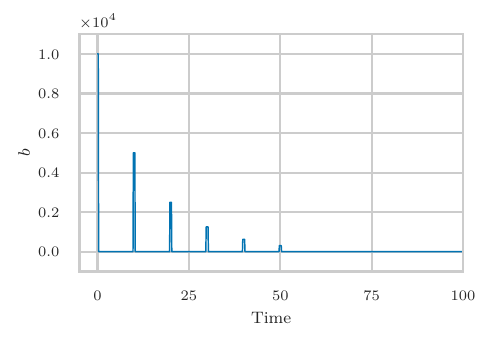}
        \caption{Temporal component $b_t$.} 
        \label{fig:dif-load1}
    \end{subfigure}
    \begin{subfigure}[htb!]{0.45\textwidth}
        \centering
        \includegraphics[trim={1.7cm 1.4cm 0.2cm 1.4cm},clip]{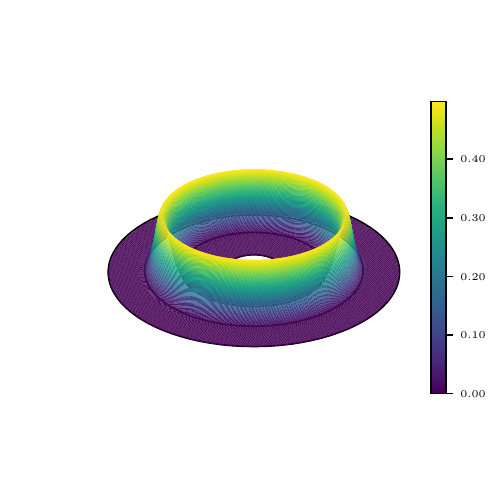}
        \caption{Spatial component $b_x$.} 
	   \label{fig:dif-load2}
    \end{subfigure}
	\caption{Diffusion problem: temporal and spatial components of the source term.}
    \label{fig:dif-load}
\end{figure}

The problem consists in solving
\begin{alignat}{2}
    \frac{\partial{u}}{\partial{t}} - \mu \Delta u &= b &\qquad \text{on} \; \Omega \times \setI \\
    u(\bm{x},t) &= 0 &\qquad \text{on} \; \Gamma \times \setI \\
    u(\bm{x},0) &= 0 &\qquad \text{on} \; \Omega. \label{eq:ht1}
\end{alignat}
We define the load using an affine decomposition of its spatial and temporal components, $b = b_x(\bm{x}) b_t(t)$. Spatially, the load is modelled as a Gaussian function centred in the middle ring of the domain, while temporally, it is defined as a sequence of point loads with varying amplitudes. In \cref{eq:heat-load} we define the load $b$, which we illustrate in \cref{fig:dif-load}.
\begin{align}
    b_x(\bm{x}) &= \cosh \left( 25 \sqrt{2 \pi} \norm{\bm{x}} \right) \exp \left( - \frac{25 \norm{\bm{x}}^2 + 50 \pi}{2} \right)\\
    b_t(t) &= \frac{A}{\sqrt{2\cdot 10^{-4}\pi}} \exp \left( {\frac{-(t-t_0)^2}{2\cdot 10^{-4}}} \right), \label{eq:heat-load}
\end{align}
with $A = 1\cdot10^4 \cdot [1,\nicefrac{1}{2},\nicefrac{1}{4},\nicefrac{1}{8},\nicefrac{1}{16},\nicefrac{1}{32}]$, and $ t_0 = 10\cdot[0,1,2,3,4,5]$. We approximate $\mathcal{L} \{b_t\}$ as the Laplace transform of a delayed impulse $\widehat{b}(s) = \e^{-t_0s}$.

The goal of this example is to compare the different possibilities of constructing the basis. We fix the sample size as $n_s = 30$ and the basis size $r=15$. We also set $\alpha=1$ and $\beta=2$. We define six different approaches to construct the reduced basis
\begin{itemize}
    \item Using $n_s$ real-valued samples $\bm{\mathsf{S}} = [\Re \{\widehat{u}_1 \}, \ldots, \Re \{\widehat{u}_{n_s} \}]$. We denote this basis as $\phi_\IR$.
    \item Using $2 n_s$ real-valued samples $\bm{\mathsf{S}} = [\Re \{\widehat{u}_1 \}, \ldots, \Re \{\widehat{u}_{n_s} \}]$. We denote this basis as $\phi_\IR^{2n_s}$.
    \item Using $n_s$ complex-valued samples $\widehat{\bm{\mathsf{S}}} = [\widehat{u}_1, \ldots, \widehat{u}_{n_s}]$, and using the complex-valued basis $\Phi$. Note that when using this basis, the solution of \cref{pr:dpr} is complex-valued. We denote this basis as $\phi_\IC$.
    \item Using $n_s$ complex-valued samples $\widehat{\bm{\mathsf{S}}} = [\widehat{u}_1, \ldots, \widehat{u}_{n_s}]$, and using the real part of the obtained complex basis $\Re \{ \Phi_\IC \}$. We denote this basis as $\phi_{\Re\IC}$.
    \item Using $n_s$ complex-valued samples arranged in a real-valued set composed by the real part and imaginary part of the frequency solution as $\bm{\mathsf{S}}^S = [\Re \{\widehat{u}_1 \}, \ldots, \Re \{\widehat{u}_{n_s}\}, \Im \{\widehat{u}_1 \}, \ldots, \Im \{\widehat{u}_{n_s}\}]$. We denote this basis as $\phi_{S}$.
    \item Using $n_s=30$ real-valued samples obtained from the standard finite element solution $\bm{\mathsf{S}} = [u_1, \ldots, u_{n_s}]$. We construct this basis following a standard reduced-basis \gls{pod} approach and denote it as $\phi_t$.
\end{itemize}

\Cref{tab:basis-types} summarizes the structural differences between the bases used in this example: their nature (real vs. complex), the number of samples used to construct them, and the snapshot set. While $\phi_t$ relies on dense sampling in the time domain, all frequency-based constructions achieve comparable or better accuracy with significantly fewer samples, and only $\phi_\IC$ requires complex-valued computations.

\begin{table}[h]
\centering
\caption{Summary of basis types used in the advection example.}
\label{tab:basis-types}
\begin{tabular}{cccc}
\toprule
\textbf{Basis Type} & \textbf{Real/Complex} & \textbf{Snapshot Set} & \textbf{Samples Used} \\
\midrule
$\phi_\IC$               & Complex & $\{\hat{u}(s_j)\}$                          & $n_s$ \\
$\phi_S$                 & Real    & $\{[\Re \hat{u}(s_j);\, \Im \hat{u}(s_j)]\}$ & $n_s$ \\
$\phi_{\Re\IC}$          & Real    & $\{\hat{u}(s_j)\}$                      & $n_s$ \\
$\phi_{\IR}$             & Real    & $\{\Re \hat{u}(s_j)\}$                      & $n_s$ \\
$\phi_{\IR}^{2n_s}$      & Real    & $\{\Re \hat{u}(s_j)\}$                      & $2n_s$ \\
$\phi_t$                 & Real    & $\{u(t_j)\}$                                & $n_s$ \\
\bottomrule
\end{tabular}
\end{table}

\begin{figure}[H]
\RawFloats
\begin{minipage}[t]{.57\textwidth}
    \includegraphics{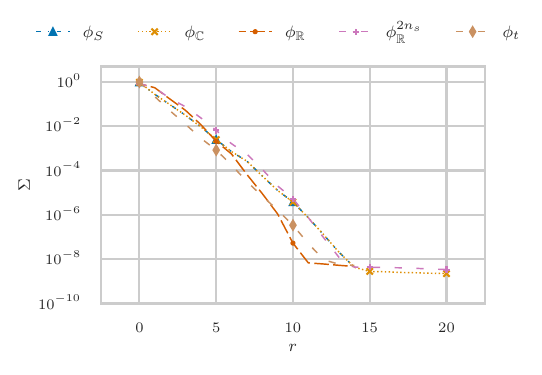}
	\caption{Diffusion problem: singular value decay for different reduced bases.} 
	\label{fig:dif-basisS}
\end{minipage}\hfill%
\begin{minipage}[t]{.385\textwidth}
    \includegraphics[width=1.11\textwidth]{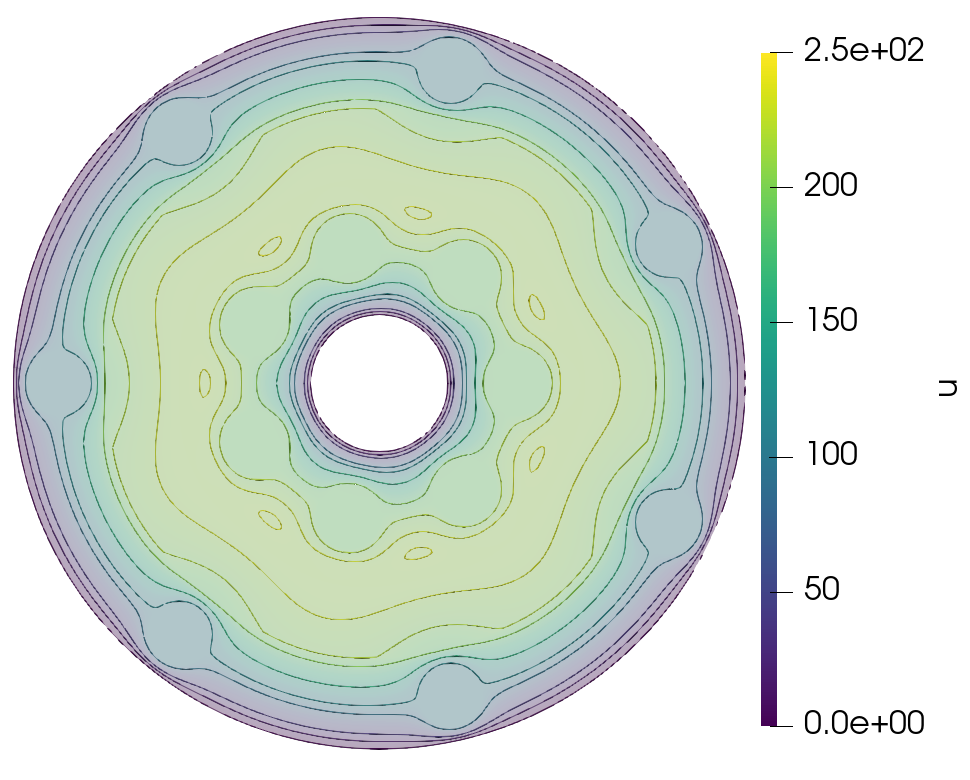}
 	\caption{Diffusion problem: reduced-order solution at \( t = 25 \).} 
	\label{fig:dif-contour}
\end{minipage}
\end{figure}

\Cref{fig:dif-contour} shows the reduced solution at $t = 25$. Qualitatively,  the frequency-based reduced-order solutions ---constructed using $\phi_\IC$, $\phi_S$, $\phi_{\Re\IC}$, and the enriched real basis $\phi_{\mathbb{R}}^{2n_s}$--- accurately reproduce the reference solution, without visible artifacts or instabilities. In contrast, the standard \gls{pod} basis $\phi_t$, built from time-domain snapshots, performs significantly worse and fails to resolve the sharp gradients characteristic of the problem. 
The singular value decay in \cref{fig:dif-basisS} further shows that $\phi_\IC$ and $\phi_S$ are nearly indistinguishable, supporting the structural equivalence discussed in \cref{ssec:complexified-basis}.

\begin{figure}[H]
	\centering
    \includegraphics{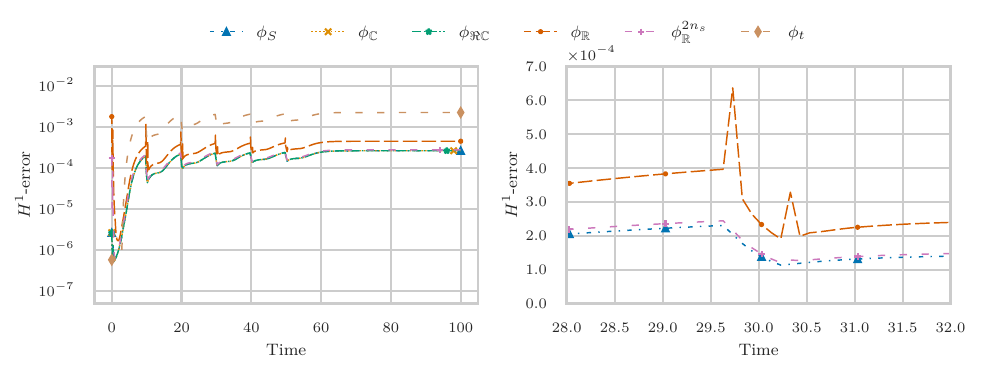}
    \caption{Diffusion problem: relative $H^1$-error over time for different basis constructions.}
	\label{fig:dif-error}
\end{figure}

\Cref{fig:dif-error} presents the evolution of the relative $H^1$-error (\cref{eq:H1-error-rel}) for all six bases. The complex-valued basis $\Phi_\IC$, the symplectic basis $\Phi_S$, and the realified basis $\Phi_{\Re\IC}$ yield nearly identical error curves, indicating that the reduced spaces they span are functionally equivalent in this setting. In contrast, the purely real-valued basis $\Phi_{\IR}$, constructed from real parts only, results in noticeably higher error unless the number of frequency samples is doubled—as done in the enriched basis $\Phi_{\IR}^{2n_s}$, which restores accuracy close to that of the complex and symplectic constructions. Finally, the time-domain basis $\Phi_t$ performs worst overall, with both slower error decay and poor resolution of transient features. These results demonstrate that preserving the analytic structure of the Laplace-transformed solutions ---either through complexification, symplectic lifting, or careful enrichment--- is essential for accurate reduced order models.

In \cref{tab:dif-time} we show the computational time for the offline and online parts of the solution including the computation of the reduced basis. The computational cost of the time-dependent finite element solution is $730.8104$ and the cost of recovering the full solution from the reduced one ($\Phi \bm{\mathsf{u}}_j^{\textrm{R}}$) for all the time steps is $58.19$.
\begin{table}[!htb]
    \centering
    \caption{Computational cost of the diffusion problem solutions.}
    {\small
    \begin{tabular}{lccccc}
        \toprule
        & $\phi_{\IR}$ & $\phi_t$ & $\phi_{\Re\IC}$ & $\phi_{\IC}$ & $\phi_{S}$ \\ \midrule  
        Frequency sampling & 254.91 & 939.81 & 254.91 & 254.91 & 254.91 \\
        Basis construction & 26.68 & 84.80 & 129.77 & 129.77 & 44.89 \\
        ROM solution & 0.49 & 0.50 & 0.78 & 0.56 & 0.49 \\
        \bottomrule
    \end{tabular}}
\label{tab:dif-time}
\end{table}

Since our goal is to obtain accurate reduced-order models at low computational cost, it is natural to consider real-valued formulations. Among these, the symplectic construction of the reduced basis $\Phi^S$ offers the best alternative: it avoids complex-valued arithmetic, it is straightforward to implement using standard real-valued linear algebra routines, and its performance comparable to the complex-valued basis.

\paragraph{Computational trade-offs.}
Both \gls{frb} and standard time-domain reduced-basis methods use Galerkin projection in the online stage, so their online complexity is essentially identical. The main difference lies in the offline stage.

For standard reduced-basis, constructing the basis requires solving the full-order time-dependent problem over $n_t$ time steps, resulting in a total cost of approximately $\mathcal{O}(n_t \cdot n_h^3)$ when using direct solvers, where $n_h$ denotes the number of spatial degrees of freedom. In contrast, the \gls{frb} method constructs the basis by solving $n_s$ independent elliptic problems in the frequency domain, each of cost $\mathcal{O}(n_h^3)$. The total offline complexity for \gls{frb} is thus $\mathcal{O}(n_s \cdot n_h^3)$.

This represents a significant reduction in offline cost when $n_s \ll n_t$, particularly since the \gls{frb} solves are fully parallelisable. Moreover, because the reduced basis is constructed independently of the time discretization, the \gls{frb} method is highly adaptable to different time integration schemes.

\subsection{Advection-dominated problem} \label{ssec:CDR}

\begin{figure}[H]
\RawFloats
\begin{minipage}[t]{.5\textwidth}
    \centering
    \includegraphics{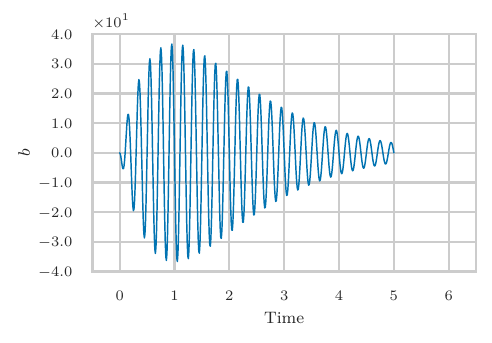}
    \caption{Advection problem: time-dependent load \( b_t(t) \).}
    \label{fig:adv-load}
\end{minipage}\hfill%
\begin{minipage}[t]{.5\textwidth}
    \centering
    \includegraphics{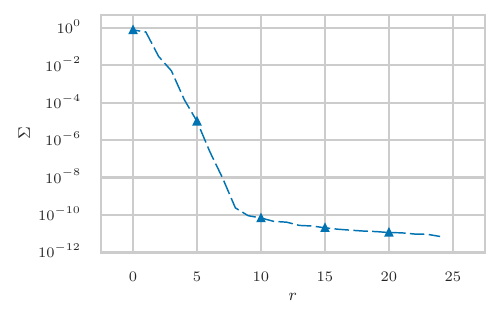}
    \caption{Advection problem: singular value decay of the reduced basis.}
    \label{fig:adv-basisS}
\end{minipage}
\end{figure}

In the last example we showcase some stabilization properties of the method. We solve the advection-diffusion-reaction problem
\begin{alignat}{2} \label{eq:advection}
    \frac{\partial{u}}{\partial{t}} + a \cdot \nabla u - \mu \Delta u + \kappa u &= b &\qquad \text{on} \; \Omega \times (0,t_{\textnormal{f}}] \\
    u(\bm{x},t) &= 0 &\qquad \text{on} \; \Gamma_D \times (0,t_{\textnormal{f}}] \\
    u(\bm{x},0) &= u_0(\bm{x};\mu) &\qquad \text{on} \; \Omega
\end{alignat}
with the $\mu=0.001$, $\kappa=0.1$, and the advection velocity $a=10-\abs{x_1}$. The domain $\Omega$ consists in a unit square with the boundary $\Gamma_D$ set over the edges parallel to the x-axis. The load and its Laplace transform are define as
\begin{align}
    b_x(\bm{x}) &= \exp \left(-\gamma^2 \left[(x_0-0.5)^2 + (x_1-0.5)^2 \right] \right), \\
    b_t(t) &= A t\e^{-t}\sin(\omega t), \quad \widehat{b}(s) = \frac{2\omega (s+1)}{((s+1)^2+\omega^2)^2}
\end{align}
with $A=-1$, $\omega=10$, and $\gamma=25$. For the frequency domain problem we set $\alpha=1$ and $\beta=2$, and we construct the basis using $30$ samples and solve the reduced problem with $n_r=15$ basis vectors. \Cref{fig:adv-load} illustrates the load $b$ and in \cref{fig:adv-basisS} we show the singular values of the reduced basis.

To assess the performance of the \gls{frb} method, we additionally compute solutions using two classical stabilization techniques: the \gls{supg} and the \gls{vms} methods. Both approaches use a stabilization term of the form $\mathcal{L}^*(v_h)\, \tau\, \mathcal{R}(u_h)$, where the stabilization parameter is defined as $\tau = \frac{h}{2 \norm{a}}$ in both cases.
In the \gls{supg} method, the residual and adjoint operators are defined as:
\begin{equation}
   \mathcal{R}(u_h) = a \cdot \nabla u_h, \qquad \mathcal{L}^*(v_h) = a \cdot \nabla v_h, 
\end{equation}
while the \gls{vms} formulation employs the full strong residual of the equation as:
\begin{equation}
    \mathcal{R}(u_h) = \frac{\partial u_h}{\partial t} + a \cdot \nabla u_h - \mu \Delta u_h + \kappa u_h - b, \qquad
    \mathcal{L}^*(v_h) = - a \cdot \nabla v_h - \mu \Delta v_h + \kappa v_h.
\end{equation}
For further details on the derivation and theoretical background of these methods, we refer the reader to \cite{Codina2000, Codina2018}.

\begin{figure}[H]
    \centering
    \includegraphics{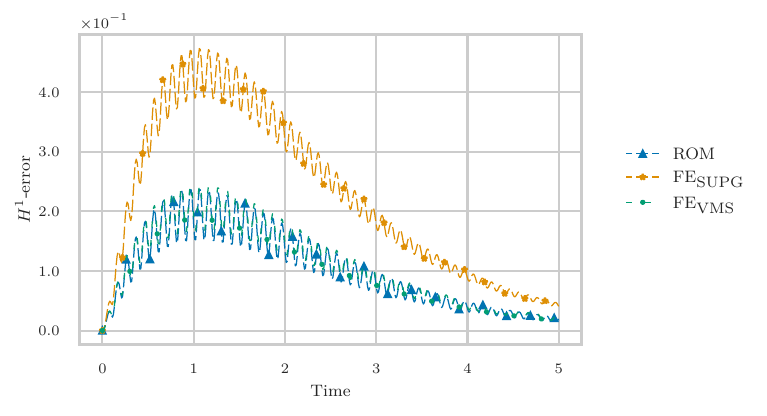}
    \caption{Advection problem: $H^1$-error over time for ROM and stabilized FEM solutions.} 
    \label{fig:adv-error}
\end{figure}

We compare the $H^1$-error (\cref{eq:H1-error}) of the three solutions ---\gls{frb}, \gls{supg} and \gls{vms}---   with respect to a high-fidelity finite element solution obtained on a finer mesh of $93356$ elements.
\Cref{fig:adv-error} displays the evolution of this error over time, showing that the \gls{frb} solution closely follows the \gls{vms} solution and achieves better accuracy than the \gls{supg} solution.

\begin{figure}[H]
	\centering
	\begin{subfigure}[htb!]{0.495\textwidth}
	    \centering
	    \includegraphics[width=1.125\textwidth]{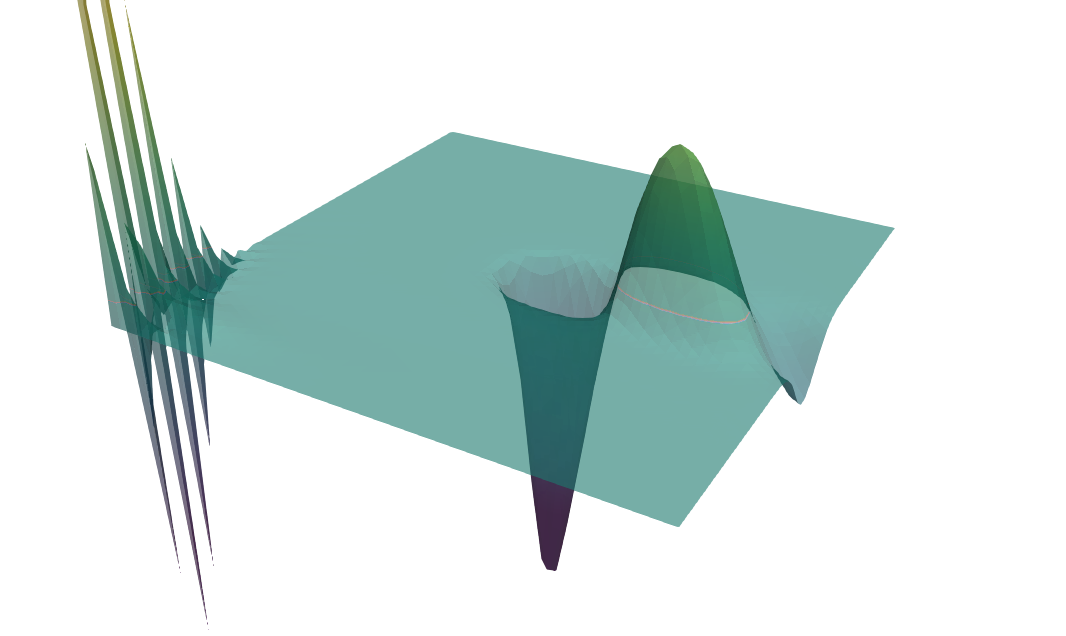}
	    \caption{Finite element solution}
    \end{subfigure}
    \begin{subfigure}[htb!]{0.495\textwidth}
	    \centering
	    \includegraphics[width=1.125\textwidth]{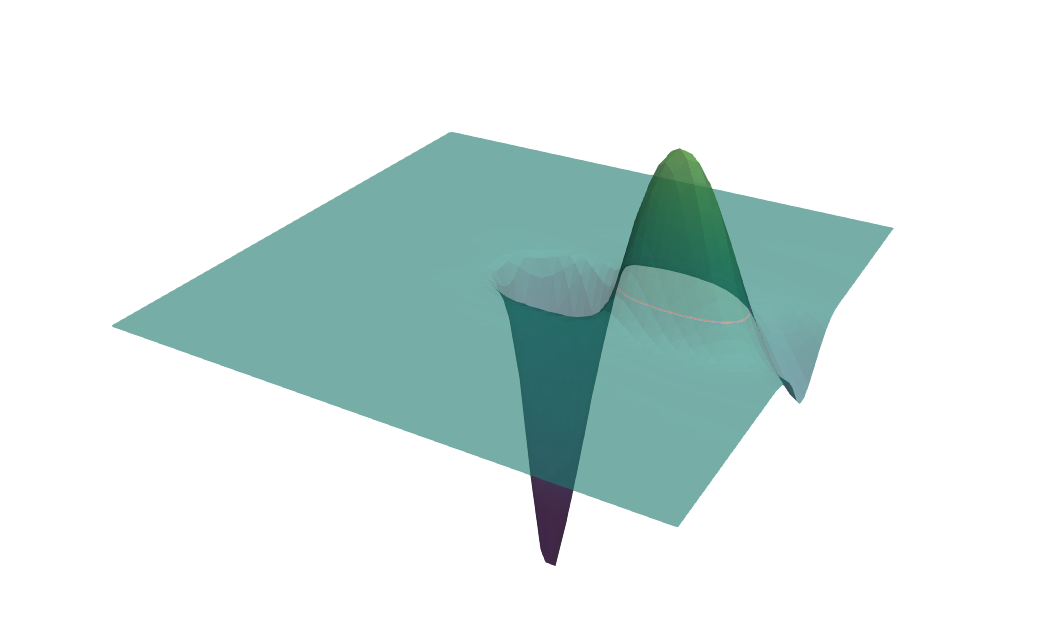}
	    \caption{VMS - Finite element solution}
    \end{subfigure} \\
	\begin{subfigure}[htb!]{0.495\textwidth}
	    \centering
	    \includegraphics[width=1.125\textwidth]{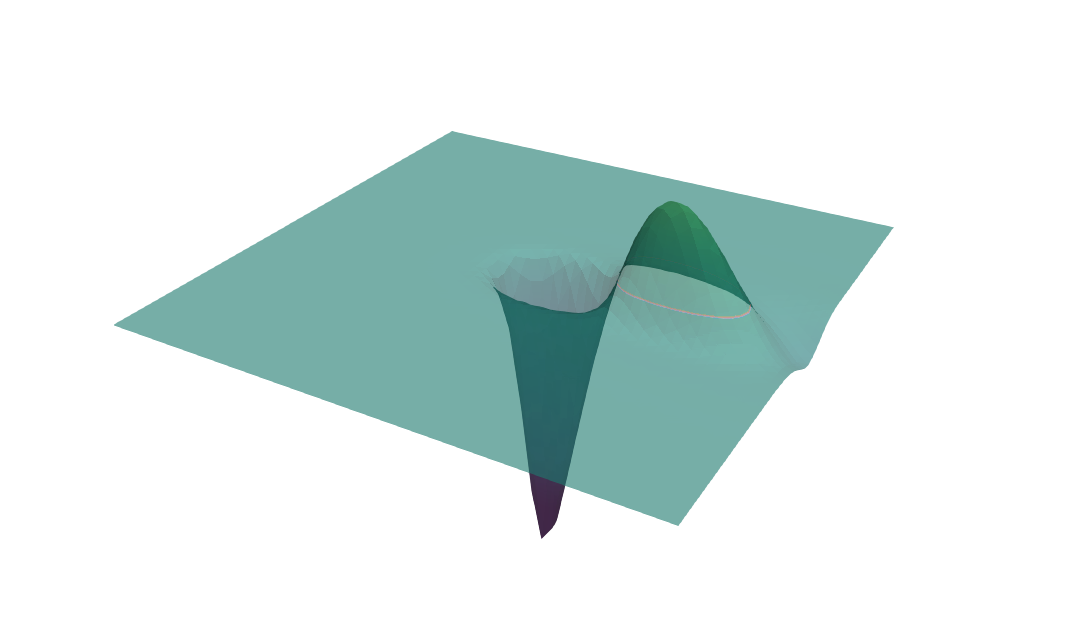}
	    \caption{SUPG - Finite element solution}
    \end{subfigure}
    \begin{subfigure}[htb!]{0.495\textwidth}
	    \centering
	    \includegraphics[width=1.125\textwidth]{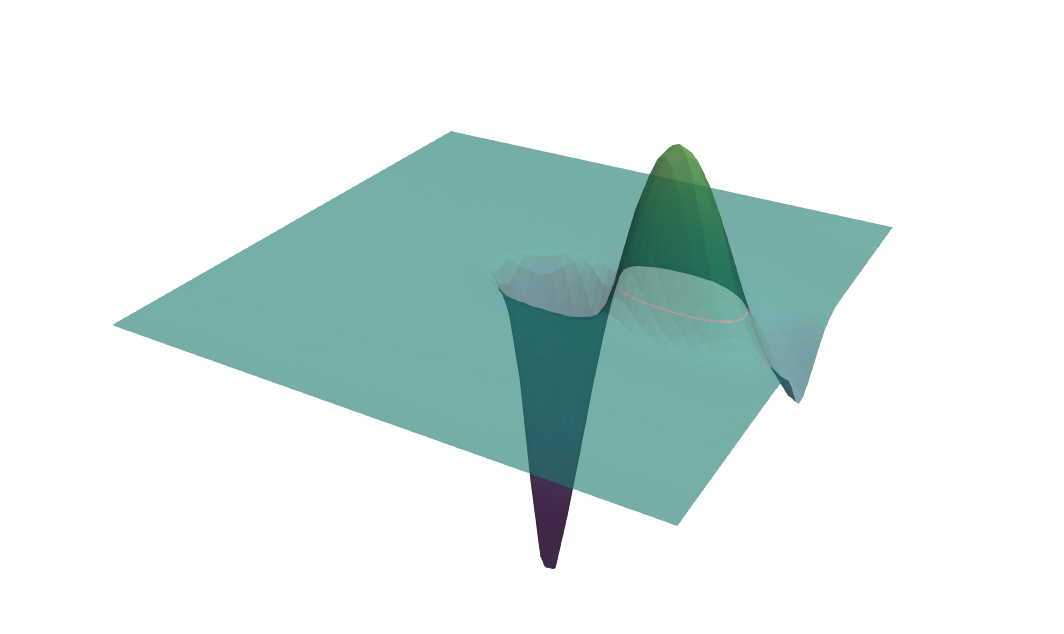}
	    \caption{ROM solution}
    \end{subfigure} \\
    \vspace{0.4cm}
    \begin{subfigure}[htb!]{\textwidth}
	    \centering
	    \includegraphics[scale=0.4]{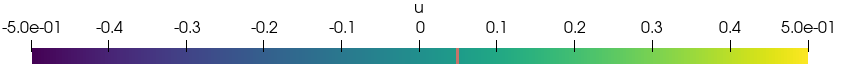}
    \end{subfigure}
    \caption{Advection problem: solution comparison at \( t = 0.33 \) for ROM, SUPG, VMS, and unstabilized FEM.}
	\label{fig:adv-comparison}
\end{figure}

In \cref{fig:adv-comparison}, we show the solutions of the finite element problems ---both stabilized and not stabilized---alongside the reduced-order model at time $t=0.66$. The classical Galerkin solution exhibits strong spurious oscillations due to the high Péclet number and the absence of stabilization. Both \gls{supg} and \gls{vms} improve stability, with the \gls{vms} method behaving less diffusive as expected (see \cite{Codina2000}). In contrast, the reduced-order solution remains smooth and free of oscillations, despite no stabilization being explicitly applied. It also proves less diffusive than the \gls{supg} solution, yielding results comparable to those of the \gls{vms} method.

This qualitative behaviour illustrates that the \gls{frb} method can suppress non-physical oscillations, likely due to the fact that it solves elliptic problems in the Laplace domain and truncates high-frequency modes through the \gls{pod}. Since the reduced basis is constructed from smooth frequency responses, unresolved temporal components ---often responsible for instabilities in convection-dominated regimes--- are naturally excluded. Unlike classical methods that rely on artificial diffusion or upwinding, the \gls{frb} achieves this stabilization intrinsically, without tuning parameters or problem-specific adjustments.

\begin{remark}

The numerical robustness of the \gls{frb} method in convection-dominated regimes is notable, as no additional stabilization is applied to the reduced-order model. This behaviour can be partly attributed to the structural properties of the Laplace-transformed formulation. By converting the original time-dependent problem into a parametric elliptic problem in the frequency domain, the method avoids direct time integration and its associated stability constraints.

In contrast to time-domain reduced-basis methods, which rely on real-valued snapshots that may contain localized transients or numerical noise, the \gls{frb} method constructs the basis from smooth, globally defined frequency-domain solutions. These solutions correspond to elliptic problems parametrized by complex frequencies, and their spatial structure often reflects integrated, steady-state-like behaviour. As a result, the reduced basis inherits a degree of regularity that improves its capacity to represent advective features without oscillation or overshoot. Moreover, the use of frequency-domain solutions inherently captures oscillatory modes and long-time behaviour, which are typically challenging to resolve using real-valued time snapshots.

While a full theoretical justification for robustness in the hyperbolic or high-Péclet regime is beyond the scope of this work, the numerical results suggest that the \gls{frb} formulation remains stable without additional stabilization mechanisms, even in cases where traditional \gls{pod}-Galerkin methods may require them.

\end{remark}
\pagebreak
\section{Summary} \label{sec:Summary}

    We propose the \textit{Frequency Reduced-Basis} method, a model reduction technique for time-dependent problems based on sampling and compression in the Laplace domain. The method follows a traditional \gls{pod} but with an important distinction: the sampling is done over the frequency domain. In the paper we provided a path for the derivation of the method following some theorems of harmonic analysis. We also show an interpretation of the constructed basis and the equivalences using some concepts of differential geometry.
    The most important features of this method are the following:
\begin{itemize}    
    \item As the method is based on the norm equivalence between time and frequency domains, the obtained reduced basis is a basis for both problems. 
	\item The solutions obtained with this reduced basis maintain good accuracy while keeping the computational cost below the high fidelity solution, even when including the sampling and construction of the basis.
    \item Since we are not modifying the \gls{pod} and the Laplace transform is a linear map, any extension of the reduced basis method proposed for stationary problems can be implemented. For example parametric reduced basis, local basis, etc.
    \item Problems with linear transformations in the source term have an equivalent reduced basis. For example, we can think of problems where the source term is time shifted, in which case the Laplace transform is invariant.
    \item Initial conditions behave as right hand side terms in the frequency domain problem, so any linear transformation on them results in equivalent reduced bases.
    \item The \gls{frb} method is compatible with any time-stepping scheme, as the reduced basis is constructed independently of the time discretization.
\end{itemize}

\paragraph{Limitations and Outlook.}
While the \gls{frb} method provides clear advantages in terms of sampling efficiency and phase accuracy, it is inherently limited to problems for which the Laplace transform yields a well-posed elliptic surrogate. Non-linear systems, or problems with strongly time-dependent coefficients or memory effects, fall outside the scope of the current framework. Extending the method to such settings remains an important direction for future work.



\bibliography{references}

\end{document}